\documentclass[12pt]{iopart}
\usepackage{amssymb,latexsym}
\usepackage{psfig}

\newcommand{\R}{\ensuremath{\mathbb{R}}}

\newcommand{\mvec}{\ensuremath{\mathbf{m}}}
\newcommand{\nvec}{\ensuremath{\mathbf{n}}}
\newcommand{\uvec}{\ensuremath{\mathbf{u}}}
\newcommand{\vvec}{\ensuremath{\mathbf{v}}}
\newcommand{\wvec}{\ensuremath{\mathbf{w}}}
\newcommand{\xvec}{\ensuremath{\mathbf{x}}}

\newcommand{\Div}{\ensuremath{\mathrm{div}}}

\newcommand{\Eh}{\ensuremath{\mathcal{E}_h}}
\newcommand{\Th}{\ensuremath{\mathcal{T}_h}}

\begin{document}

\title[Simulation of Combustion in Porous Media]
      {Hybrid Discretization Methods\\ for Transient Numerical Simulation\\
       of Combustion in Porous Media}

\author{Peter Knabner\footnote[7]{knabner@am.uni-erlangen.de} and Gerhard Summ}

\address{Institute of Applied Mathematics, Martensstra{\ss}e 3, D-91058 Erlangen, Germany}

\begin{abstract}
We present an algorithm for the numerical solution of the equations governing
combustion in porous inert media.
The discretization of the flow problem is performed 
by the mixed finite element method,
the transport problems are discretized by a cell-centered finite volume method.
The resulting nonlinear equations are lineararized with Newton's method, 
the linearized systems are solved with a multigrid algorithm.
Both subsystems are recoupled again in a Picard iteration.
Numerical simulations based on a simplified model show
how regions with different porosity stabilize the reaction zone 
inside the porous burner
\end{abstract}


\submitto{\CTM}


%
\section{Introduction}
\label{sec:intro}
%
In recent years, the request for low-emission combustion systems
has led to the development of a new burner concept:
combustion in porous media.
In these burners a premixed gas-air mixture is constrained to flow through 
and combust within the pores of a porous medium, typically a ceramic foam.
Due to the presence of the solid matrix, the temperature in the combustion zone
can be controlled such that the emission of $\mathrm{NO}_x$ and other 
pollutants can be lowered.
Further advantages of this technique are mentioned in \cite{Brenner/ea}.
A list of technical applications can be found in \cite{Trimis/Durst}.

Several authors performed numerical simulations to support the development
of porous burners.
Most of these simulations are based on mathematical models
that are spatially one-dimensional \cite{Hanamura/Echigo} or
stationary in time \cite{Brenner/ea},\cite{Mohamad/Ramadhyani/Viskanta},
or one-dimensional and stationary \cite{Hsu/Howell/Matthews}.
The employed numerical algorithms are often simple extensions 
of existing codes designed for other applications
and use finite difference or finite volume methods for the discretization 
of the governing equations.
Some of these simulations consider the pressure distribution and the velocity
field of the flow as given. 
In \cite{Brenner/ea} and \cite{Mohamad/Ramadhyani/Viskanta}
the pressure distribution and the velocity field are computed 
from the continuity equation and a momentum equation.
Then pressure and velocity are coupled using the SIMPLE-algorithm.
We will present a discretization of the two-dimensional transient problem,
which includes the computation of the pressure distribution 
and the velocity field with the mixed finite element method.
By this means, we obtain an accurate approximation of the velocity field,
which influences the convective transport of heat and chemical species.

In the next section, we introduce the governing equations, 
initial and boundary conditions of the problem.
We restrict our considerations here to a simple irreversible one-step
reaction mechanism.
An extension of the numerical algorithm to more realistic reaction systems 
is directly possible.
In Section \ref{sec:disc} we discuss the discretization of these equations
and algorithms used for the solution of the resulting algebraic equations.
Finally, in Section \ref{sec:ex} we present numerical calculations based on the
simplified model, which indicate the stabilizing effect of regions
with different porosity on the location of the combustion zone.

%
\section{Mathematical model}
\label{sec:model}
%
\mathindent6em%
The simulation of the flow and transport processes inside the porous medium 
based on models that resolve the complex pore structure in detail 
would be computationally too expensive.
Hence these simulations are based on homogenized or volume-averaged equations,
where the fluid and the solid are considered as a pseudohomogeneous medium.
We restrict our considerations here to the simplest possible model,
which can be used to describe the flow of the gas in the porous burner 
and the transport of heat and chemical species.
This model is sufficiently realistic to present the essential features
of our algorithm.

%
\subsection{Model equations}
The flow in the porous medium is governed by the Darcy--Forchheimer equation
\[
\frac{\mu}{k} \, \uvec + \beta_{\mathrm{Fo}} \, \rho \, |\uvec| \, \uvec
+ \nabla p = 0 \: ,
\]
where $| \cdot |$ denotes the Euclidean norm, the continuity equation
\[
\phi \, \frac{\partial \rho}{\partial t} 
+ \Div(\rho \uvec) = 0 
\]
and the ideal gas law 
\[
\rho = p \, \frac{W}{R_0 \, T} \: .
\]
The unknowns here are the pressure $p$, the density $\rho$ 
and the volumetric flow rate $\uvec$ of the gas.
Porosity $\phi$, permeability $k$ 
and Forchheimer coefficient $\beta_{\mathrm{Fo}}$ of the porous medium,
viscosity $\mu$, molecular weight $W$ and temperature $T$ of the gas 
mixture, and the universal gas constant $R_0$ are given.
Assuming $\rho>0$ and introducing new variables $S = |p| p$ and 
$\mvec = |\rho| \uvec$ these equations can be transformed into
\begin{eqnarray}				\label{Da-Fo.eq}
\left( \alpha + \beta |\mvec| \right) \mvec + \nabla S = 0 \: , \\
\phi \, \partial_t \rho(S) + \Div(\mvec) = 0 \: , \label{cont.eq}
\end{eqnarray}
where 
\[ \gamma := \frac{W}{R_0 \, T} \; , \quad
   \alpha := \frac{2 \, \mu}{\gamma \, k} \; , \quad
   \beta := \frac{2 \, \beta_{\mathrm{Fo}}}{\gamma} \; , \]
and the equation of state $\rho=\rho(S)$ is defined by
\[ 
\rho(S) := \gamma \frac{S}{\sqrt{|S|}} \: .
\] 

The solid matrix of the porous burner typically consists of ceramic foams
having high heat transfer rates. Furthermore the specific surface of these
foams is very high. Therefore we can assume thermal equilibrium between
the gas and the solid. The temperature of this pseudo-homogeneous medium
is governed by the (effective) heat equation
\mathindent3em%
\begin{equation}				\label{heat-eff}
\phi c_p \partial_t ( \rho T ) + (1-\phi) c_s \partial_t ( \rho_s T )
+ \Div \big( c_p \mvec \, T - \lambda_{\mathrm{eff}} \nabla T \big)
= \phi Q \dot{r} + (1-\phi) F_Q \; ,
\end{equation}
\mathindent6em%
where the effective heat conductivity
$\lambda_{\mathrm{eff}} := \phi \lambda_g + (1-\phi) \lambda_s$
is a mean value of the heat conductivities $\lambda_g$ of the gas 
and $\lambda_s$ of the solid.
More realistic representations of $\lambda_{\mathrm{eff}}$, 
including effects of radiation and dispersion (see \cite{Brenner/ea}), 
could also be used.
Furthermore $c_p$ denotes the heat capacity for constant pressure of the gas,
$c_s$ the specific heat capacity and $\rho_s$ the density of the solid.
$Q \dot{r}$ is the heat produced by the combustion reaction and
$F_Q$ the power density of an external heat source needed for ignition.

Consider the simplest case only,
we restrict the reaction model to an irreversible one-step reaction mechanism
(e.g.\ methane oxidation 
 $\mathrm{CH}_4 + 2 \mathrm{O}_2 \longrightarrow 
  \mathrm{CO}_2 + 2 \mathrm{H}_2 \mathrm{O}$).
Then it suffices to consider the conservation equation of the reactant $R$
\begin{equation}			\label{species-reactant}
\phi \partial_t ( \rho y ) +
\Div \! \left( \mvec y - \phi D \nabla y \right) = - \phi \dot{r} \; ,
\end{equation}
where $y$ denotes the mass fraction of the reactant and the sink term $-\dot{r}$
models the consumption of the reactant by the chemical reaction. 
Following the Arrhenius model $\dot{r}$ is given by
\[
\dot{r} = B \rho y \exp \left( -\frac{E}{R_0 T} \right)
\]
with the frequency factor $B$ and the activation energy $E$.

%
\subsection{Initial and boundary conditions}
Since (\ref{Da-Fo.eq}) and (\ref{cont.eq}) are equivalent to a single
parabolic equation, we have to prescribe initial conditions for $S$
and boundary conditions for $S$ or $\mvec$ only.
Hence we have to provide the following initial conditions for the problem 
consisting of equations (\ref{Da-Fo.eq})--(\ref{species-reactant}):
\[
  S(\cdot,0) = S_0 \; \mbox{ in} ~ \Omega \; , \quad
  T(\cdot,0) = T_0 \; \mbox{ in} ~ \Omega \; , \quad
  y(\cdot,0) = y_0 \; \mbox{ in} ~ \Omega \; ,
\]
where $\Omega \subset \R^2$ is the domain corresponding to the porous burner.
Then the initial condition for the mass flux $\mvec(\cdot,0)$ 
can be calculated from (\ref{Da-Fo.eq}).

\begin{table}[ht]%
\caption{Boundary conditions}%
\label{tab:bc:comb}%
$\begin{array}{lrcl}
& & & \\[-1.5ex]
\mbox{Flow problem}: & & & \\[0.5ex]
\mbox{Inflow boundary } \Gamma_{\mathrm{I}}: & \mvec \cdot \nvec & = & m_b \\
\mbox{Outflow boundary } \Gamma_{\mathrm{O}}: & S & = & S_b \\
\mbox{Burner wall } \Gamma_{\mathrm{W}}: & \mvec \cdot \nvec & = & 0 \\
\mbox{Line of symmetry } \Gamma_{\mathrm{S}}: &
\mvec \cdot \nvec & = & 0 \\[1ex]
\mbox{Heat equation}: & & & \\[0.5ex]
\mbox{Inflow boundary } \Gamma_{\mathrm{I}}: &
\lambda_{\mathrm{eff}} \nabla T \cdot \nvec
& = & c_p \left( T - T_b \right) m_b \\
\mbox{Outflow boundary } \Gamma_{\mathrm{O}}: & \nabla T \cdot \nvec & = & 0 \\
\mbox{Burner wall } \Gamma_{\mathrm{W}}: &
\lambda_{\mathrm{eff}} \nabla T \cdot \nvec & = & h_b \left( T_b-T \right) \\
\mbox{Line of symmetry } \Gamma_{\mathrm{S}}: &
\nabla T \cdot \nvec & = & 0 \\[1ex]
\mbox{Reactant conservation equation}: \quad & & & \\[0.5ex]
\mbox{Inflow boundary } \Gamma_{\mathrm{I}}: & y & = & y_b \\
\mbox{Outflow boundary } \Gamma_{\mathrm{O}}: & \nabla y \cdot \nvec & = & 0 \\
\mbox{Burner wall } \Gamma_{\mathrm{W}}: & \nabla y \cdot \nvec & = & 0 \\
\mbox{Line of symmetry } \Gamma_{\mathrm{S}}: &
\nabla y \cdot \nvec & = & 0
\end{array}$
\end{table}

The boundary conditions are listed in Table~\ref{tab:bc:comb}.
Note that the mixed boundary conditions for the heat equation model
thermal equilibrium at the inflow boundary and Newton's law of cooling
at the burner wall. 
The coefficients $h_b$ and $T_b$ denote the heat transfer coefficient
and the ambient temperature, resp.

%
\section{Numerical solution algorithm}
\label{sec:disc}
%
Using the Rothe method, we first discretize equations 
(\ref{Da-Fo.eq})--(\ref{species-reactant}) in time with the 
implicit Euler method.
Hence, in each time step we have to solve the following coupled system:
\begin{eqnarray}				\label{semid.DaFo.eq}
\left( \alpha + \beta |\mvec| \right) \mvec + \nabla S = 0 \: , \\
						\label{semid.cont.eq}
\frac{\phi}{\Delta t} (\rho(S) - \rho^-) + \Div(\mvec) = 0 \: , \\
						\label{semid.heat-eff}
\frac{\phi c_p}{\Delta t} ( \rho(S) T - \rho^- T^- )
+ \frac{(1-\phi)c_s}{\Delta t} (\rho_s T - \rho_s^- T^-) \\[0.5ex] \nonumber
\hspace{9.1em} + \: \Div \left( c_p \mvec T - \lambda_{\mathrm{eff}} \nabla T \right)
= \phi Q \dot{r} + (1-\phi) F_Q \: , \\	\label{semid.species-reactant}
\frac{\phi}{\Delta t} ( \rho(S) y - \rho^- y^- )
+ \Div \left( \mvec y - \phi D \nabla y \right) 
= - \phi \dot{r} \: .
\end{eqnarray}
Here $\rho^-$, $T^-$ and $y^-$ denote the values of the unknowns 
from the last time step.

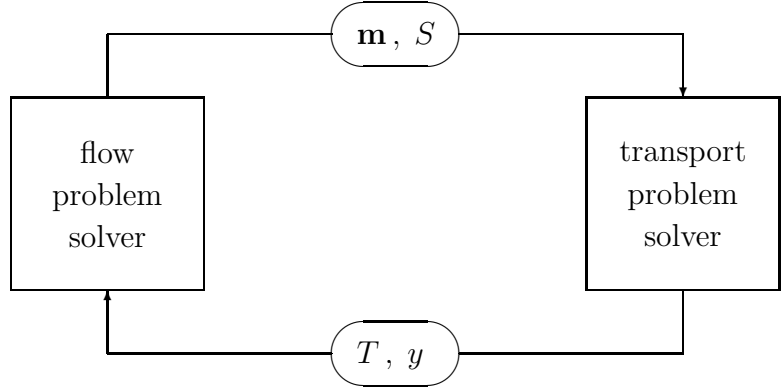
\begin{figure}[ht]
\centerline{
\setlength{\unitlength}{0.85mm}%
\begin{picture}(120,60)%
\put(0,0){\parbox{0mm}{\begin{center} \end{center}}}%
\put(0,15){\framebox(30,30){%
		\parbox{30mm}{%
			\begin{center}%
                        	flow\\problem\\solver%
			\end{center}}}}
\put(15,45){\line(0,1){10}}
\put(15,55){\line(1,0){35}}
\put(60,55){\oval(20,10)}
\put(54,53.5){$\mvec \, , \; S$}
\put(70,55){\line(1,0){35}}
\put(105,55){\vector(0,-1){10}}
\put(90,15){\framebox(30,30){%
		\parbox{30mm}{%
			\begin{center}%
                        	transport\\problem\\solver%
			\end{center}}}}
\put(105,15){\line(0,-1){10}}
\put(105,5){\line(-1,0){35}}
\put(60,5){\oval(20,10)}
\put(54,3.5){$T \, , \; y$}
\put(50,5){\line(-1,0){35}}
\put(15,5){\vector(0,1){10}}
\end{picture}}
\caption[]{Decoupling of flow and transport problem}
\label{Picard}
\end{figure}

The coupling between the flow problem 
(governed by equations (\ref{semid.DaFo.eq}) and (\ref{semid.cont.eq}))
and the transport problems (governed by 
equations (\ref{semid.heat-eff}) and (\ref{semid.species-reactant}))
is rather weak. Therefore we can decouple the flow problem 
and the transport problems in each time step and compute their solutions
$(\mvec,S)$ and $(T,y)$ alternatingly in a Picard-iteration,
as sketched in \Fref{Picard}.
This strategy allows us to use different methods for the discretization
of the flow and the transport problem.

%
\subsection{Solution of the flow problem}
The spatial discretization of both problems is based on a partition
${\mathcal{T}}_h$ of the domain $\Omega$ into triangular elements $K$.
We denote by $\mathcal{E}_h$ the set of edges of $\mathcal{T}_h$.
Respecting the partition of $\partial \Omega$ into 
Dirichlet boundary $\Gamma_D$ (where $S = S_b$ is given) and 
Neumann boundary $\Gamma_N$ (where $\mvec \cdot \nvec = m_b$ is given),
$\mathcal{E}_h$ can be subdivided into three disjoint subsets
$\mathcal{E}_h = \mathcal{E}_h^I \cup \mathcal{E}_h^D \cup \mathcal{E}_h^N$.
Here $\mathcal{E}_h^I := \left\{ e \in \Eh \bigm| 
				 e \not \subset \partial \Omega \right\}$ 
is the set of inner edges, 
$\mathcal{E}_h^D := \left\{ e \in \Eh \bigm| e \subset \Gamma_D \right\}$ 
is the set of edges lying on the Dirichlet boundary and
$\mathcal{E}_h^N := \left\{ e \in \Eh \bigm| e \subset \Gamma_N \right\}$ 
is the set of edges lying on the Neumann boundary.
For every $K \in {\mathcal{T}}_h$ we denote by $P_k(K)$, $k \ge 0$, 
the set of polynomials of degree $\le k$ on $K$.
Analogously, we define $P_k(e)$ for $k \ge 0$ and $e \in \mathcal{E}_h$.

To preserve the form given by (\ref{semid.DaFo.eq}) and (\ref{semid.cont.eq}),
and to get a good approximation for the mass flux density $\mvec$, 
which appears also in equations (\ref{semid.heat-eff}) and (\ref{semid.species-reactant}),
we use the mixed finite element method on Raviart--Thomas elements of 
lowest order (see e.g.\ \cite{Brezzi/Fortin}) for the spatial 
discretization of (\ref{semid.DaFo.eq}) and (\ref{semid.cont.eq}).
Thus $\mvec$ is approximated by $\mvec_h \in V_h$, where
\mathindent3em
\[
  V_h := RT_0(\Omega;{\mathcal{T}}_h) :=
  \left\{ \vvec_h \in H(\mathrm{div};\Omega) \bigm| 
          \vvec_h|_K \in RT_0(K) \quad 
          \mbox{for all } K \in {\mathcal{T}}_h \right\}
\]
\mathindent6em
and $RT_0(K)$ is defined by
\[
  RT_0(K) := \big( P_0(K) \big)^2 + {x \choose y} P_0(K) \; ,
\]
and $S$ is approximated by $S_h \in Q_h$, where
\[
  Q_h := \left\{ q_h \in L^2(\Omega) \bigm| 
                 q_h|_K \in P_0(K) \quad 
                 \mbox{for all } \, K \in {\mathcal{T}}_h \right\} \; .
\]
Note that any $\vvec_h|_K \in RT_0(K)$ is uniquely defined by the 
flux across the edges of $K$
$$
  v_e := \int_e \vvec_h|_K \cdot \nvec_e \, ds \quad , 
  \quad e \subset \partial K \; ,
$$
where $\nvec_e$ is an arbitrarily oriented unit normal vector 
to $e \in \mathcal{E}_h$.
The property $\vvec_h \in H(\mathrm{div};\Omega)$ 
requires the continuity of these flux values.
Using the basis 
$\left\{ \wvec_e \right\}_{e \in {\mathcal{E}}_h} $ of $V_h$ satisfying
$$
  \int_f \wvec_e \cdot \nvec_f \, ds = 
  \left\{ \begin{array}{ccc} 
		1 \: , & \mathrm{if} & e=f \: , \\
                0 \: , & \mathrm{if} & e \neq f \: , \end{array} \right. 
$$
the corresponding degrees of freedom are given by the fluxes 
$v_e$, $e \in \mathcal{E}_h$.

To take into account the flux boundary conditions on $\Gamma_N$,
we consider the subspace $V_h^{m_b,\Gamma_N}$ of $V_h$,
which is defined by
\[ V_h^{m_b,\Gamma_N} := \left\{ \vvec_h \in V_h \Bigm| 
   \int_e \vvec_h \cdot \nvec \, ds = \int_e m_b \, ds \quad
   \mbox{for all } e \in \mathcal{E}_h^N \right\} \: . \]
If $m_b \equiv 0$ on $\Gamma_N$, we obtain the space $V_h^{0,\Gamma_N}$.

Then the discrete mixed formulation reads as follows:\\
Find $(\mvec_h,S_h) \in V_h^{m_b,\Gamma_N} \times Q_h$ such that
for every $(\vvec_h,q_h) \in V_h^{0,\Gamma_N} \times Q_h$
\mathindent3em
\begin{eqnarray}			\label{disc_mixed:DaFo}
\int_\Omega \left( \alpha + \beta |\mvec_h| \right)
	    \left( \mvec_h \cdot \vvec_h \right) \, dx
- \int_\Omega \Div (\vvec_h) \, S_h \, dx 
+ \int_{\Gamma_D} S_b \left( \vvec_h \cdot \nvec \right) ds = 0 \: , \\ 
					\label{disc_mixed:cont}
\int_\Omega \frac{\phi}{\Delta t} \rho(S_h) \, q_h \, dx 
+ \int_\Omega \Div \left( \mvec_h \right) q_h \, dx 
= \int_\Omega \frac{\phi}{\Delta t} \rho^- \, q_h \, dx \; .
\end{eqnarray}
\mathindent6em

Unfortunately, the systems of algebraic equations resulting from the mixed
finite element discretization are difficult to solve numerically.
Therefore we apply the following implementational technique,
called hybridization to our discrete mixed formulation:
We eliminate the continuity constraints in the definition of $V_h$
and enforce the required continuity instead through 
additional equations involving Lagrange multipliers
defined on the edges $e \in \mathcal{E}_h$.
Thus we replace $V_h$ by
\mathindent3em
\[
  W_h := RT_{-1}(\Omega;{\mathcal{T}}_h) :=
  \left\{ \vvec_h \in \left(L^2(\Omega)\right)^2 \bigm| 
          \vvec_h|_K \in RT_0(K) \quad 
          \mbox{for all } K \in {\mathcal{T}}_h \right\} \, ,
\]
and $V_h^{m_b,\Gamma_N}$ by the corresponding subspace $W_h^{m_b,\Gamma_N}$
of $W_h$.
In addition, we define the space of Lagrange multipliers by
\[ 
   \Lambda_h^{S_b,\Gamma_D} := 
   \left\{ \lambda_h \in L^2 (E_h) \Bigm|
   \lambda_h|_e \in P_0(e) ~ \forall e \in \mathcal{E}_h , 
	   \int_e (\lambda_h-S_b) \, ds = 0 
		  ~ \forall \, e \in \mathcal{E}_h^D \right\} \; ,
\]
where $E_h=\cup_{e \in \mathcal{E}_h} e$.
Then the hybridized mixed formulation reads as:

Find $(\mvec_h,S_h,\mu_h) 
      \in W_h^{m_b,\Gamma_N} \times Q_h \times \Lambda_h^{S_b,\Gamma_D}$, 
such that for every
$(\vvec_h, q_h, \lambda_h) 
 \in W_h^{0,\Gamma_N} \times Q_h \times \Lambda_h^{0,\partial \Omega}$ 
\mathindent0em
\begin{eqnarray}			\label{hybr_mixed:DaFo}	
\int_\Omega \! \left( \alpha + \beta |\mvec_h| \right)
	    \left( \mvec_h \cdot \vvec_h \right) dx
- \sum_{K \in \mathcal{T}_h}
      \int_K \! \Div \left( \vvec_h \right) S_h \, dx 
+ \sum_{K \in {\mathcal{T}}_h} \int_{\partial K} \! \mu_h 
			\left( \vvec_h \cdot \nvec_K \right) ds = 0 \: , \\
					\label{hybr_mixed:cont}
\int_\Omega \frac{\phi}{\Delta t} \rho(S_h) \, q_h \, dx 
+ \sum_{K \in \mathcal{T}_h} \, 
    \int_K \Div \left( \mvec_h \right) q_h \, dx
= \int_\Omega \frac{\phi}{\Delta t} \rho^- \, q_h \, dx \: , \\
					\label{hybr_mixed:Lagr}
\sum_{K \in \mathcal{T}_h} \, \int_{\partial K} 
	\lambda_h \left( \mvec_h \cdot \nvec_K \right) ds = 0 \: .
\end{eqnarray}
\mathindent6em%

The solutions $\mvec_h$ and $S_h$ of 
(\ref{hybr_mixed:DaFo})--(\ref{hybr_mixed:Lagr})
coincide with the solutions $\mvec_h$ and $S_h$ of 
(\ref{disc_mixed:DaFo})--(\ref{disc_mixed:cont}) (cf.~\cite{Knabner/Summ}).
Therefore we are allowed to use the same notation for them.
We note that the additionally computed Lagrange multipliers can be used 
to construct a more accurate approximation for $S$
(see e.g.\ \cite{Brezzi/Fortin} and the references cited there).

Note that every $\vvec_h \in W_h$ is uniquely defined by the 
degrees of freedom $v_{K,e}$, which are given by
\[
  v_{K,e} = \int_e \vvec_h|_K \cdot \nvec_K \, ds ~,
  \quad K \in \mathcal{T}_h , ~ e \subset \partial K \: ,
\]
where $\nvec_K$ denotes the unit outer normal of $K$.
The corresponding basis vectors are denoted by $\wvec_{K,e}$.
Thus the unknown functions $\mvec_h$, $S_h$ and $\mu_h$ can be represented by
\mathindent1em%
\[
\mvec_h = \sum_{K \in \mathcal{T}_h} 
    	  \sum_{e \in \partial K} m_{K,e} \, \wvec_{K,e} \quad , \quad
S_h = \sum_{K \in \mathcal{T}_h} S_K \, \chi_K \quad \mathrm{and} \quad
\mu_h = \sum_{e \in \mathcal{E}_h} \mu_e \, \chi_e \: ,
\]
\mathindent6em%
where $\chi_K$ and $\chi_e$ denote the characteristic functions
of an element $K$ and an edge $e$, resp..

Employing the basis functions 
$\wvec_{K,\bar{e}}$ for $\bar{e} \notin \mathcal{E}_h^N$ and 
$K \in \mathcal{T}_h$ with $\bar{e} \subset \partial K$,
$\chi_K$ for $K \in \mathcal{T}_h$, and
$\chi_e$ for $e \in \mathcal{E}_h^I$ as test functions, 
we obtain from (\ref{hybr_mixed:DaFo})--(\ref{hybr_mixed:Lagr})
the following system of nonlinear equations.
\mathindent0em%
\begin{eqnarray}
F_{K,\bar{e}} := \!		\label{F_(K,e)}
\int_K \! \left( \alpha + \beta |\mvec_h| \right)
       \left( \mvec_h \cdot \wvec_{K,\bar{e}} \right) d\xvec
 - S_K + \mu_{\bar{e}} = 0 \: , 
& \bar{e} \notin \mathcal{E}_h^N , 
  K \in \mathcal{T}_h, \bar{e} \subset \partial K \: , \\
F_K := 				\label{F_K}
\frac{\int_K \phi\, \gamma \, d\xvec}{\Delta t} \, \frac{S_K}{\sqrt{|S_K|}}
+ \sum_{e \subset \partial K} m_{K,e}
- \frac{\int_K \phi\, \rho^- \, d\xvec}{\Delta t} = 0 \: , \:
& K \in \mathcal{T}_h \: , \\	\label{F_e}
F_{\bar{e}} := m_{K,\bar{e}} + m_{K',\bar{e}} = 0 \: ,
& e \in \mathcal{E}_h^I \: .
\end{eqnarray}
\mathindent6em%
For every $K \in \mathcal{T}_h$ the basis functions $\wvec_{K,e}$
($e \subset \partial K$) of $W_h$ vanish in $\Omega \setminus K$.
Therefore the degrees of freedom $m_{K,e}$ for $e \subset \partial K$ affect 
only $\mvec_h|_K = \sum_{e \in \partial K} m_{K,e} \, \wvec_{K,e} =: \mvec_K$.
Hence, for every $K \in \mathcal{T}_h$, the degrees of freedom $m_{K,e}$ 
for $e \subset \partial K$ and $S_K$ appear only in equations
(\ref{F_(K,e)}) and (\ref{F_K}) belonging to the element $K$.
This fact can be exploited to compensate the main drawback of hybridization,
the introduction of additional degrees of freedom:
For every $K \in \mathcal{T}_h$, we can eliminate 
the approximate mass flux values $m_{K,e}$ 
and the approximate element values $S_K$ from the global system
by solving locally the systems consisting of (\ref{F_(K,e)}) and (\ref{F_K}).

Owing to the nonlinearity of (\ref{F_(K,e)}) and (\ref{F_K}) 
it is not possible to find a closed form solution 
for $m_{K,e}$ ($e \subset \partial K$) and $S_K$.
Nevertheless, using monotonicity arguments it can be shown that 
the local subsystems are uniquely solvable and their Jacobians are invertible
(see \cite[Section~4.3]{Knabner/Summ}).
Thus it is possible to compute $m_{K,e}$ ($e \subset \partial K$) and $S_K$
by means of Newton's method during the assembling procedure.
The elimination of $m_{K,e}$ ($e \subset \partial K$) and $S_K$ introduces 
nonlinearity into the originally linear equations (\ref{F_e}).
Using Newton's method to solve these equations we have to
compute the Jacobian of the remaining global system.
By the implicit function theorem this requires again the invertibility
of the Jacobians of the local subsystems.

The linearized global systems are solved by a multigrid method.
Chen \cite{Chen} showed that -- at least for linear elliptic problems -- 
the equations after elimination of flux and element variables correspond 
to the equations resulting from a nonconforming discretization 
using the Crouzeix--Raviart ansatz space.
Therefore we can employ the intergrid transfer operators developed 
for the Crouzeix--Raviart ansatz space \cite{Braess/Verfuerth} 
to construct a multigrid algorithm.

%
\subsection{Solution of the transport problem}
Using the equation of state $\rho=\rho(S)$, we need the values of
the temperature $T$ and auxiliary variable $S$ to compute the density $\rho$. 
Therefore we choose the same ansatz space for $T_h$ as for $S_h$, 
i.e., $T_h := \sum_{K \in \Th} T_K \chi_K \in Q_h$. 
Furthermore $\rho$, $T$ and $y$ appear in the reaction rate $\dot{r}$.
Hence we approximate the mass fraction $y$ of the reactant by 
$y_h := \sum_{K \in \Th} y_K \chi_K \in Q_h$, too.
Since the fluxes corresponding to $T$ or $y$ do not enter further equations,
we do not need to compute them explicitly using the mixed finite element method.
Instead, we employ a cell-centered finite volume scheme for the spatial 
discretization of (\ref{semid.heat-eff}) and (\ref{semid.species-reactant}).
Note that this finite volume scheme can be derived 
from the mixed finite element method (see \cite{Baranger/Maitre/Oudin}).

Integrating (\ref{semid.heat-eff}) and (\ref{semid.species-reactant}) 
over a cell $K \in \Th$, applying the divergence theorem and 
replacing the continuous unknowns $T$ and $y$ by their discrete approximations
$T_h$ and $y_h$ yields the following equations for all $K \in \Th$:
\mathindent2em%
\begin{eqnarray}	
\frac{\int_K \!\phi c_p \, d\xvec}{\Delta t} \,
	 ( \rho_K T_K - \rho_K^- T_K^- )	\label{volld.heat-eff}
  + \frac{\int_K \!(1-\phi) c_s \rho_s d\xvec}{\Delta t} \, T_K
  - \frac{\int_K \!(1-\phi) c_s \rho_s^- d\xvec}{\Delta t} \, T_K^- \\
\hspace{5em} + \sum_{e \subset \partial K} 	\nonumber
    \int_e \left( c_p \mvec_h \, T_h 
    - \lambda_{\mathrm{eff}} \nabla T_h \right) \cdot \nvec_K \, d\sigma
= \int_K \! \phi Q \dot{r} \, d\xvec 	
 + \int_K \! (1-\phi) F_Q \, d\xvec \: , \\
\frac{\int_K \!\phi \,d\xvec}{\Delta t} ( \rho_K y_K - \rho_K^- y_K^-)
+ \sum_{e \subset \partial K}			\label{volld.species-reactant}
 \int_e \left( \mvec_h y_h - \phi D \nabla y_h \right) \cdot \nvec_K \, d\sigma
= - \int_K \! \phi \dot{r} \, d\xvec \: .
\end{eqnarray}
\mathindent6em%

To obtain algebraic equations, we have to provide quadrature formulas for
the integrals appearing in (\ref{volld.heat-eff}) and (\ref{volld.species-reactant}).
The terms $\int_K \phi \, d\xvec$ and $\int_K (1-\phi) F_Q \, d\xvec$
are evaluated by means of the midpoint rule, i.e.,
$$ \int_K \phi \, d\xvec \approx |K| \, \phi(b_K) \quad \mbox{and} \quad
   \int_K (1-\phi) F_Q \, d\xvec \approx |K| \, (1-\phi(b_K)) F_Q(b_K) \: , $$
where $b_K$ is the center of gravity of $K$.
The remaining integrals over $K$ depend on one or several unknowns.
Of course, we use the corresponding element value for these unknowns.
Taking into account the Arrhenius law and the temperature dependence of
coefficient functions $c_p$, $c_s$ and $\rho_s$, 
we obtain the following approximations for these integrals:
\mathindent1em%
\begin{eqnarray*}
\int_K \phi \dot{r} \, d\xvec =
\int_K \phi B \rho_h y_h \exp \left( -\frac{E}{R_0 T_h} \right) \, d\xvec
\approx |K| \, \phi(b_K) \, B \, \rho_K \, y_K 
            \exp \! \left( -\frac{E}{R_0 T_K} \right) \: , \\
\int_K \phi c_p \, d\xvec \approx |K| \, \phi(b_K) \, c_p(T_K) \, , \\
\int_K (1-\phi) c_s \rho_s^{(-)} \, d\xvec \approx
|K| \, (1-\phi(b_K)) \, c_s(T_K) \, \rho_s(T_K^{(-)}) \: .
\end{eqnarray*}
\mathindent6em%

The approximation of integrals over the edges $e \in \Eh$ is much more 
difficult. Using the generalizing notation $z_h \in \{y_h,T_h\}$,
we have to find approximations of
$$ \int_e \left( c_z \mvec_h z_h - D_z \nabla z_h \right) \cdot \nvec_K \, d\sigma 
   \; , $$
where $z_h \in Q_h$ is piecewise constant.
In particular, $z_h$ is not continuous across interior edges. 
Consequently, $z_e := z_h|_e$ is not uniquely defined and 
$\nabla z_h|_e$ is not defined at all.
For an element $K \in \Th$ and an interior edge $e \in \Eh^I$ let
$\mathrm{Nb}(K,e) \in \Th$ be the element, which shares the edge $e$ with $K$.

We start with the consideration of interior edges $e \in \mathcal{E}_h^I$.
In order to discretize the convective flux in stable way, we use 
an upwind scheme (see e.g.\ \cite[Section~7]{Eymard/Gallouet/Herbin}), i.e.,
$$ \int_e c_z \left( \mvec_h \cdot \nvec_K \right) z_h \, d\sigma 
   \approx m_{K,e}^+ \, c_z|_K \, z_K 
         + m_{K,e}^- \, c_z|_{\mathrm{Nb}(K,e)} \, z_{\mathrm{Nb}(K,e)} \: . $$
Here $m_{K,e}^+ := \max \{ m_{K,e} , 0 \}$
 and $m_{K,e}^- := \min \{ m_{K,e} , 0 \}$.

The diffusive flux is approximated by some form of difference quotient.
Since the elements $K$ may be of different size and shape
and the value of the diffusion coefficient $D_z$ may vary 
from element to element, some points have to be respected
to obtain a stable discretization.
Similar to the computation of the effective heat conductivity 
of layered materials, we have to employ some kind of geometric mean 
for the diffusion coefficient. 
The derivation of the cell-centered finite volume method from the 
mixed finite element method yields the following discretization
(cf.\ \cite[Section 9]{Eymard/Gallouet/Herbin}):
\[ 
\int_e D_z \, \nabla z_h \cdot \nvec_K \, d\sigma \approx
|e| \frac{z_{\,\mathrm{Nb}(K,e)}-z_K}{d_e(D_z)} \; ,
\]
where
\[
d_e(D_z) = \left( \frac{d_{K,e}}{D_z|_K} 
         + \frac{d_{\mathrm{Nb}(K,e),e}}{D_z|_{\mathrm{Nb}(K,e)}} \right) \; .
\]
and $d_{K,e}$ denotes the distance between the center of the circumcircle 
of $K$ and the edge $e$.
Since the center of the circumcircle of a triangle with an obtuse angle
lies outside of the triangle, such obtuse-angled triangles have to be avoided
in the triangulation $\mathcal{T}_h$.

Finally, we describe the approximation of the fluxes across boundary edges.
In the case of Dirichlet boundary conditions ($z=z_b$ on $e$) 
we employ the approximation
\[
\int_e \left( c_z \mvec_h z_h - D_z \nabla z_h \right) \cdot \nvec_K \, d\sigma
\approx c_z(z_K) m_{K,e} z_{b,e} - |e| \, D_z(b_K) \frac{z_{b,e}-z_K}{d_{K,e}}
\: 
\]
Here $z_{b,e}$ is the mean value of $z_b$ on $e$, approximated by 
$$ z_{b,e} := \frac{1}{|e|} \int_e z_b \, d\sigma 
      \approx \left( z_b(n_1(e))+z_b(n_2(e)) \right) / 2 \: , $$
where $n_i(e)$, $i=1,2$, denote the vertices of the edge $e$.

In the case of Neumann boundary conditions $\nabla z \cdot \nvec = q_b$ 
is given. To evaluate the convective flux, we choose $z_e \approx z_K$.
Hence we obtain
\mathindent0.5em%
\begin{eqnarray*}
\int_e \left( c_z \mvec_h z_h - D_z \nabla z_h \right) \cdot \nvec_K \, d\sigma
\approx c_z(z_K) m_{K,e} z_K - D_z(b_K) \int_e q_b \, d\sigma \\
\hspace{13em} \approx c_z(z_K) m_{K,e} z_K 
            - |e| \, D_z(b_K)  \left( q_b(n_1(e))+q_b(n_2(e)) \right) / 2 \: . 
\end{eqnarray*}
\mathindent6em%

Finally, we consider the case of mixed boundary conditions, given in the form
\[
z + \frac{1}{\sigma} (D_z \nabla z \cdot \nvec ) = g \quad \iff
D_z \nabla z \cdot \nvec = \sigma (g-z) \: .
\]
In this case we employ the approximation
\[
\int_e \left( c_z \mvec_h z_h - D_z \nabla z_h \right) \cdot \nvec_K \, d\sigma
\approx c_z(z_K) m_{K,e} z_K - |e| \, \sigma_e (g_e - z_K) \; .
\]
Here $\sigma_e$ and $g_e$ are defined like $z_{b,e}$ above.
Note that the cell-centered finite volume method can be extended 
to more general elliptic operators, including discontinuous 
matrix diffusion coefficients (see \cite[Section 11]{Eymard/Gallouet/Herbin}).

After approximating the integrals in (\ref{volld.heat-eff}) and 
(\ref{volld.species-reactant}) as above, we arrive at a coupled system 
of nonlinear equations.
Again, we propose to use Newton's method for the linearization of these equations
and to use a multigrid method for the solution of the linearized equations.
For the computations presented below, we used the trivial injection operator 
$I_k^t$ and its transpose as intergrid transfer operators
and obtained satisfactory convergence rates.
As there may be cases where the $V$-cycle does not converge for these trivial
intergrid transfer operators, we mention also the weighted interpolation 
operator $I_k^w$, which has been proposed in \cite{Kwak/Kwon/Lee}.

%
\section{Stabilization of the reaction zone}
\label{sec:ex}
\begin{figure}[ht]
\centerline{%
\psfig{figure=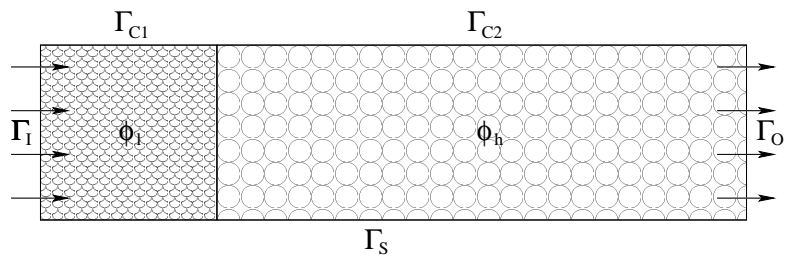,width=0.8\textwidth}}
\caption[]{Prototype of a porous burner}
\label{prototype}
\end{figure}

%
In order to demonstrate the facilities of the discretization methods 
presented above, we implemented them in the framework of the software toolbox 
\textsf{ug} \cite{Bastian/ea}.
By means of this algorithm, we studied, how regions with varying porosity can
help to control the position of the reaction zone inside the porous burner.
To this end we consider the simplified model sketched in Figure \ref{prototype}
of a porous burner prototype developed by Trimis and Durst \cite{Trimis/Durst}.
Experiments showed that the reaction zone stabilizes at the interface between
the preheating region with smaller pores (or lower porosity $\phi_\mathrm{l}$)
and the combustion region with larger pores or (or higher porosity $\phi_\mathrm{h}$).
It was suggested that this stabilization effect appears when the pore diameter
in the preheating region is lower than a quenching diameter observed 
in tube experiments.
But the pore diameter is a microscopic parameter, which does not appear in the 
macroscopic (averaged) equations (\ref{Da-Fo.eq})--(\ref{species-reactant}).
Therefore the numerical simulations presented in \cite{Brenner/ea} implement 
this stabilization effect by permitting chemical reactions only in the 
combustion region.
Our simulations show, that this stabilization effect can be explained
with the aid of macroscopic parameters only.
Of course, these simulations can not yield quantitatively realistic results,
since they are based on the simplified model consisting of equations 
(\ref{Da-Fo.eq})--(\ref{species-reactant}).
Nevertheless we expect to obtain qualitatively correct results concerning
basic features of porous burners like the stabilization effect of regions 
with varying porosity.

%
\subsection{Description of the problem}
The values of the coefficients used in our simulations are listed in
Table \ref{tab:coeff}.
For the Forchheimer parameter $\beta_{\mathrm{Fo}}$, we use Ergun's relation
$\beta_{\mathrm{Fo}} = c_{\mathrm{F}} / \sqrt{k}$, where 
$c_{\mathrm{F}}$ is a constant and $k$ the permeability of the porous medium.

\begin{table}[ht]%
\caption{Values of coefficient functions}%
\label{tab:coeff}%
\centerline{\begin{tabular}{llllll}
& & & & & \\
\multicolumn{2}{l}{Gas mixture} & 
\multicolumn{2}{l}{Solid} & 
\multicolumn{2}{l}{Reaction rate} \\[1ex] 
$\mu$ & $3.18 \cdot 10^{-5}~\mathrm{Pa~s}$ &
$c_{\mathrm{F}}$ & $0.55$ & $R_0$ & $8.314~\mathrm{J~/~mol~K}$ \\
$D$ & $8.2 \cdot 10^{-5}~\mathrm{kg~/~m~s}$ \hspace*{2em} & 
$\rho_s$ & $3970~\mathrm{kg~/}~\mathrm{m}^3$ \hspace*{2em} & 
$B$ & $1.8 \cdot 10^8~1~/~\mathrm{s}$ \\
$c_p$ & $1005~\mathrm{J~/~kg~K}$ & $c_s$ & $765~\mathrm{J~/~kg~K}$ & 
$E$ & $125600~\mathrm{J~/~mol}$ \\
$\lambda_g$ & $0.049~\mathrm{W~/~K~m}$ & $\lambda_s$ & $36~\mathrm{W~/~K~m}$ &
$Q$ & $5.0 \cdot 10^7~\mathrm{J~/~kg}$ \\
$W$ & $0.028~\mathrm{kg~/~mol}$ & & & 
\end{tabular}}
\end{table}

In order to study the influence of varying porosity (and permeability) 
on the stability and position of the reaction zone,
we compare the results of numerical simulations for the following choices 
of $\phi$ and $k$:
\begin{enumerate} \renewcommand{\labelenumi}{\alph{enumi})}
\item Different values of $\phi_{\mathrm{l}}$ and $\phi_{\mathrm{h}}$:
    \begin{equation}					\label{diff.por}
    \begin{array}{rcl}
    \phi(\xvec) & = & \left\{ \begin{array}{ll}
		  	0.3 & \hspace{3.85em} \mbox{ for } x_1 < 0.08 \: , \\
		  	0.8 & \hspace{3.85em} \mbox{ for } x_1 > 0.08 \\
	              \end{array} \right. \quad \mbox{and} \\[3ex]
    k(\xvec) & = & \left\{ \! \begin{array}{ll}
		  	1.0 \cdot 10^{-8} \, \mathrm{m}^2 
			& \mbox{ for } x_1 < 0.08 \: , \\
		  	1.0 \cdot 10^{-7} \, \mathrm{m}^2
			& \mbox{ for } x_1 > 0.08 \: . \\
	   	   \end{array} \right. \hspace{-1em}
    \end{array}
    \end{equation}
\item Constant low value $\phi \equiv \phi_{\mathrm{l}}$:
    \begin{equation}					\label{low.por}
    \phi(\xvec) = 0.3 \quad \mbox{for all } \xvec \quad \mbox {and } \quad
    k(\xvec) = 1.0 \cdot 10^{-8} \, \mathrm{m}^2 \quad \mbox{for all } \xvec \: .
    \end{equation}
\item Constant high value $\phi \equiv \phi_{\mathrm{h}}$:
    \begin{equation}					\label{high.por}
    \phi(\xvec) = 0.8 \quad \mbox{for all } \xvec \quad \mbox{and } \quad
    k(\xvec) = 1.0 \cdot 10^{-7} \, \mathrm{m}^2 \quad \mbox{for all } \xvec \: .
    \end{equation}
\end{enumerate}

We start the simulation with constant initial values
\[ S_0 \equiv 10266755625~\mathrm{Pa}^2 ~ , \quad 
   y_0 \equiv 0.0 \quad \mbox{and} \quad T_0 \equiv 298~\mathrm{K} \: . \]
Changing the flux boundary condition
$\mvec(\xvec,t) \cdot \nvec = m_b(\xvec,t)$ 
at the inflow boundary $\Gamma_{\mathrm{I}}$ linearly
from $m_b(\cdot,t) = 0~\mathrm{kg~/~m~s}$ for $t<50$
to $m_b(\cdot,t) = -0.2~\mathrm{kg~/~m~s}$ for $t>60$
a flow field is generated.
The homogeneous flux boundary conditions at the burner wall 
$\Gamma_{\mathrm{W}} = \Gamma_{\mathrm{C}1} \cup \Gamma_{\mathrm{C}2}$ 
and at the symmetry line $\Gamma_{\mathrm{S}}$, 
and the Dirichlet boundary condition $S(\xvec,t)=S_b(\xvec,t) \equiv  S_0$
are kept fixed.
At the same time ($t \in [50,60]$) the  Dirichlet boundary condition
for the mass fraction of the reactant at the inflow boundary
is increased linearly from $y_b=0$ to $y_b=0.05$.
In order to start the reaction a point heat source (modelling a glow igniter)
of power density $100\,000 \: \mathrm{W\:/}\:\mathrm{m}^3$ is located at
$\xvec = (0.1,0.07)$ for $t < 150$.
The ambient temperature $T_b$, which appears in the mixed boundary conditions
for the heat equation is kept fixed at $T_b = 298~\mathrm{K}$.
For $t>150$ the heat transfer coefficient at the burner wall 
$\Gamma_{\mathrm{W}}$ is set to
\mathindent4em%
\[ h_b = 1500~\mathrm{W~/~K~m}^2 \quad \mbox{at }~ \Gamma_{\mathrm{C}1} 
   \quad \mbox{ and } \quad
   h_b = 100~\mathrm{W~/~K~m}^2 \quad \mbox{at }~ \Gamma_{\mathrm{C}2} \: . \]
\mathindent6em%
After $t=150$ the boundary conditions are kept fixed,
the power density of the heat source is $0$.
The simulation is continued, until a numerically steady state is reached. 
This state depends strongly on the choice for the porosity and permeability 
values.

%
\subsection{Comparison of results}
\begin{figure}[ht]
\begin{center}
\begin{tabular}{cc}
\psfig{figure=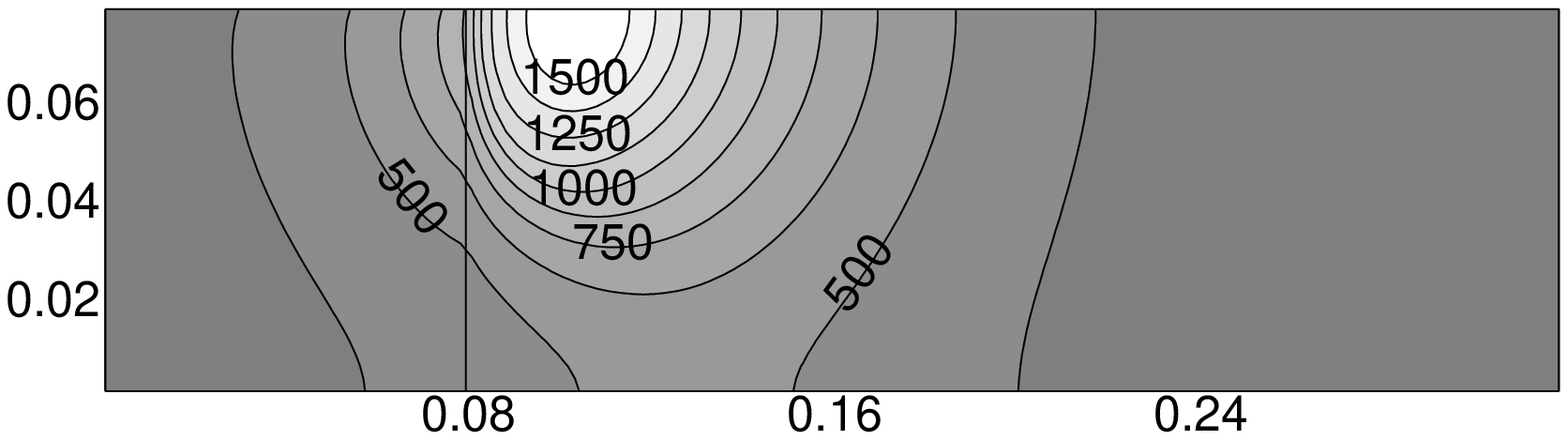,width=0.46\textwidth} &
\psfig{figure=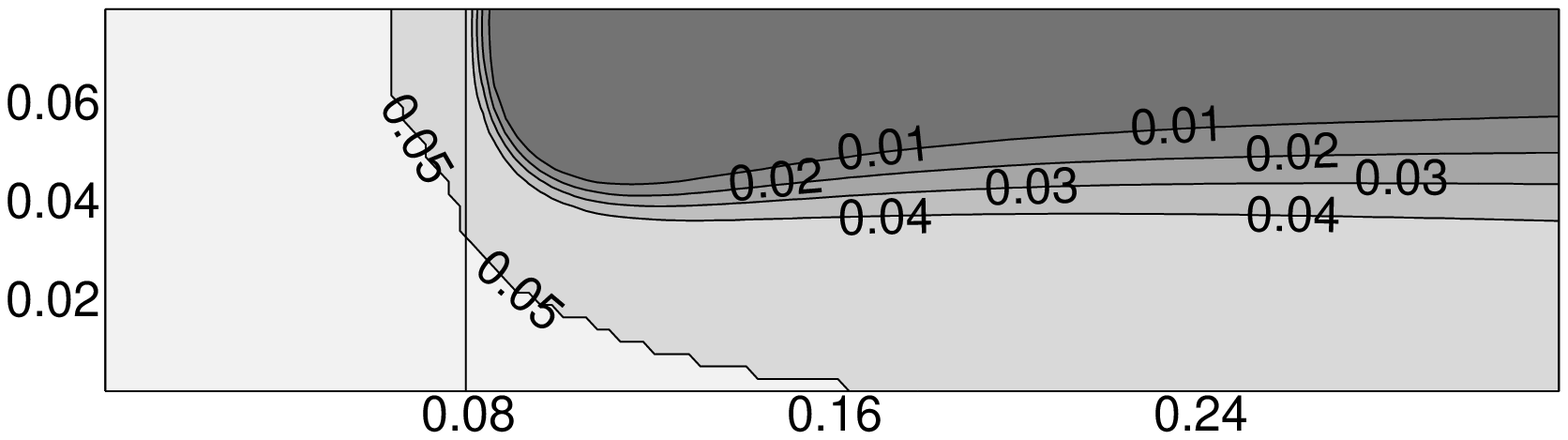,width=0.46\textwidth} \\[1ex]
\psfig{figure=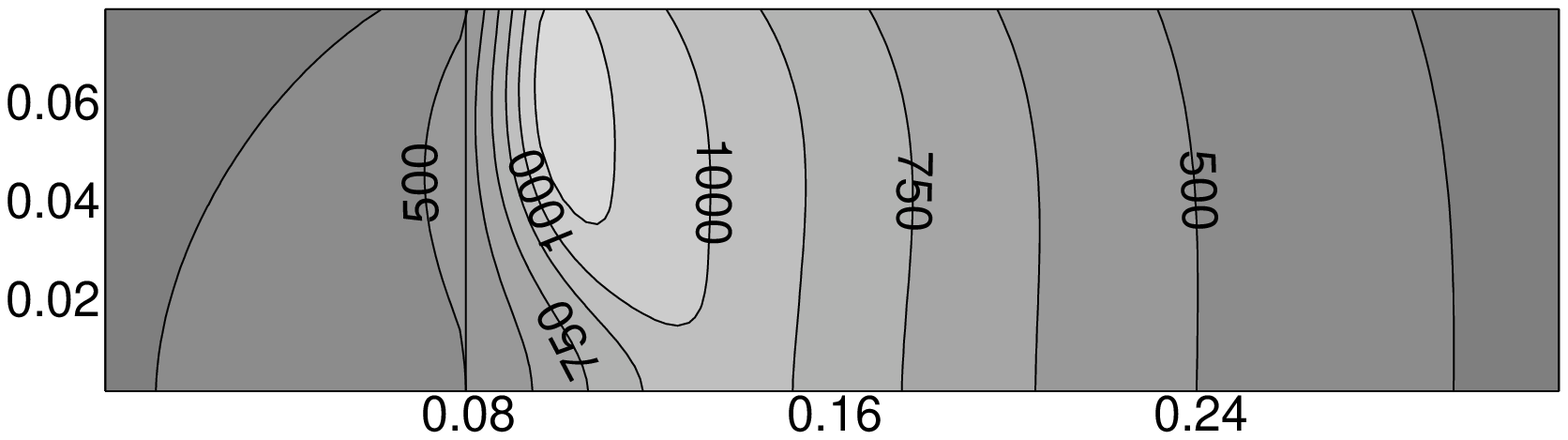,width=0.46\textwidth} &
\psfig{figure=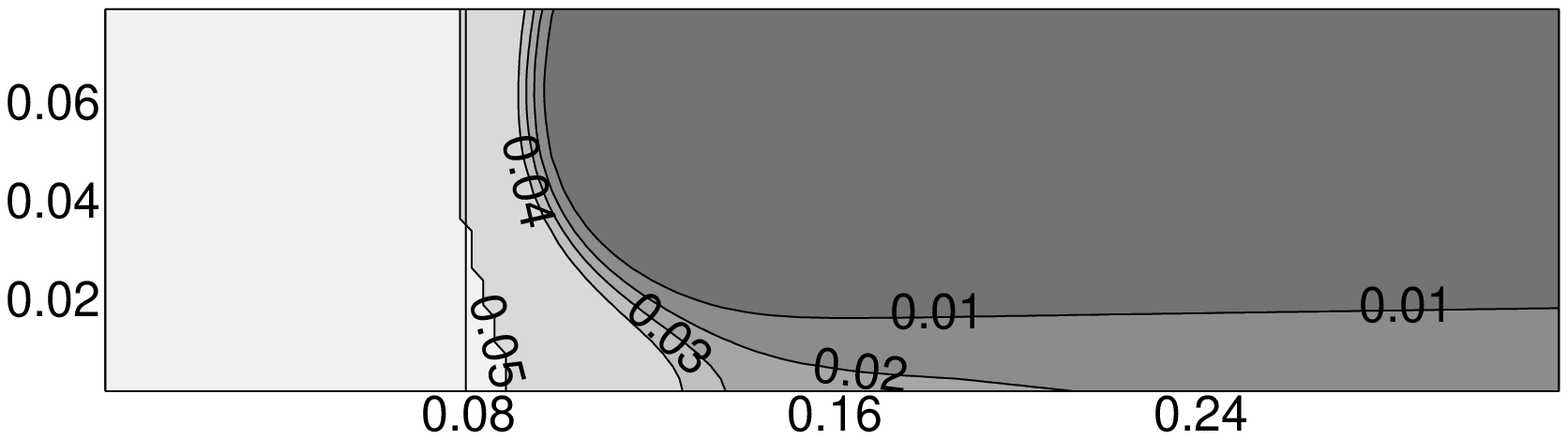,width=0.46\textwidth} \\[1ex]
\psfig{figure=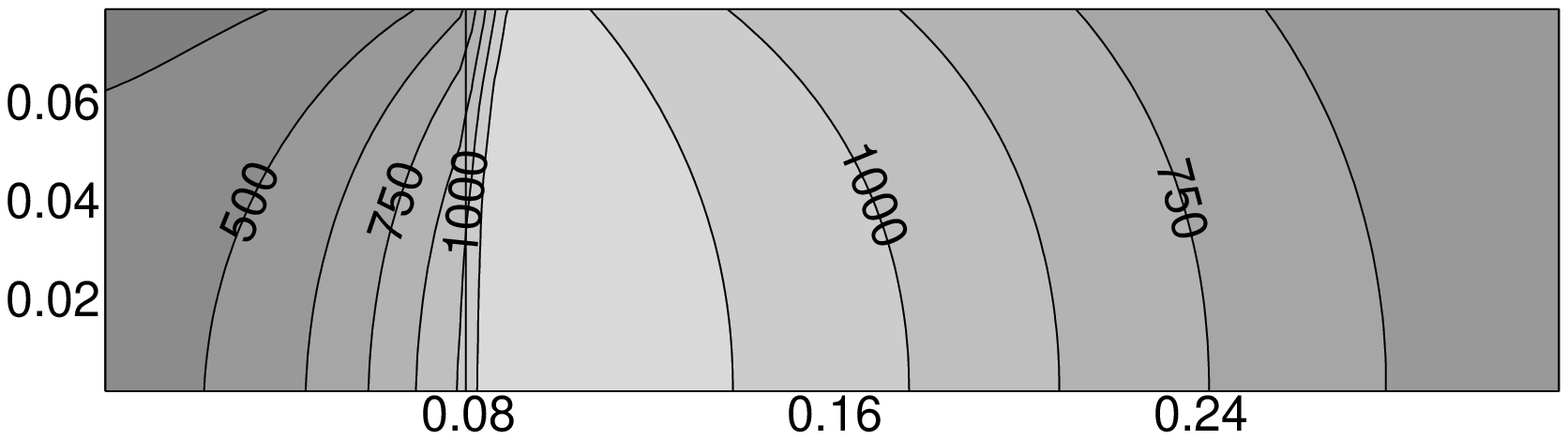,width=0.46\textwidth} &
\psfig{figure=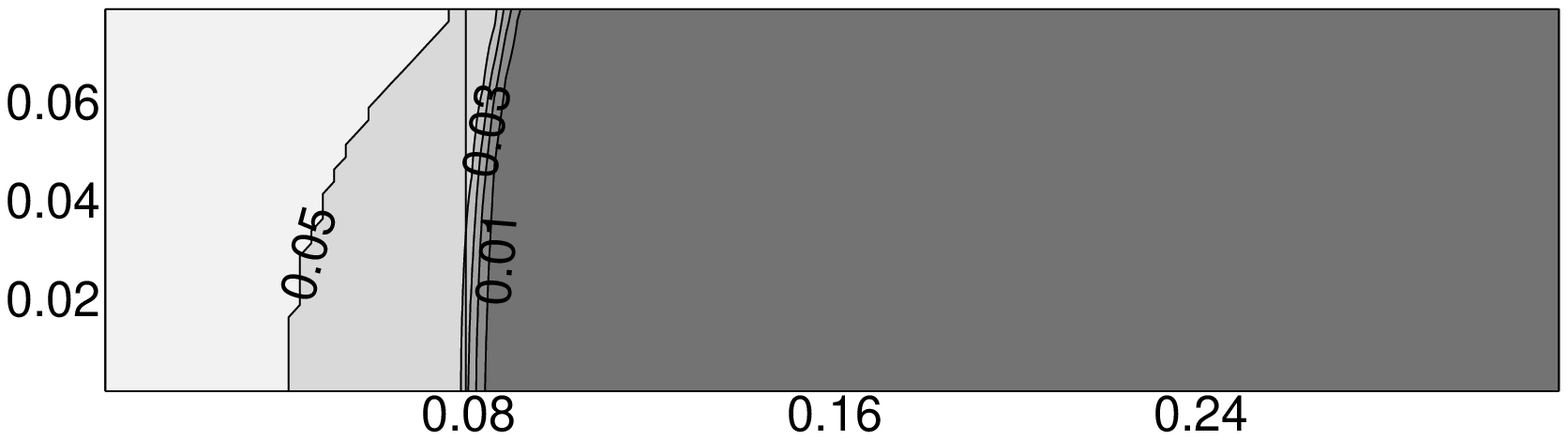,width=0.46\textwidth} \\[1ex]
\psfig{figure=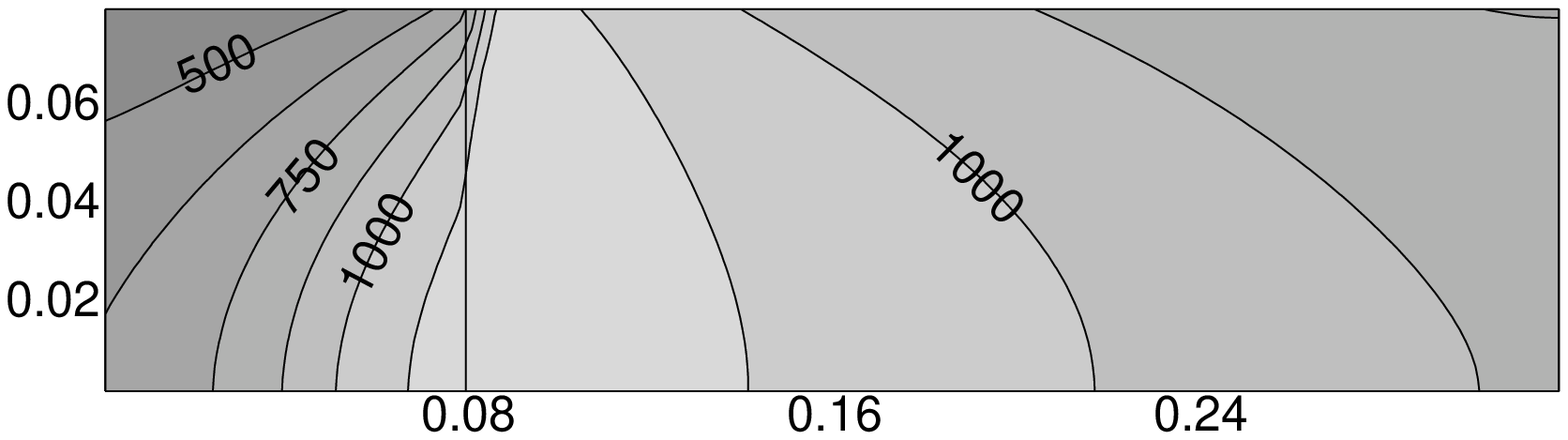,width=0.46\textwidth} &
\psfig{figure=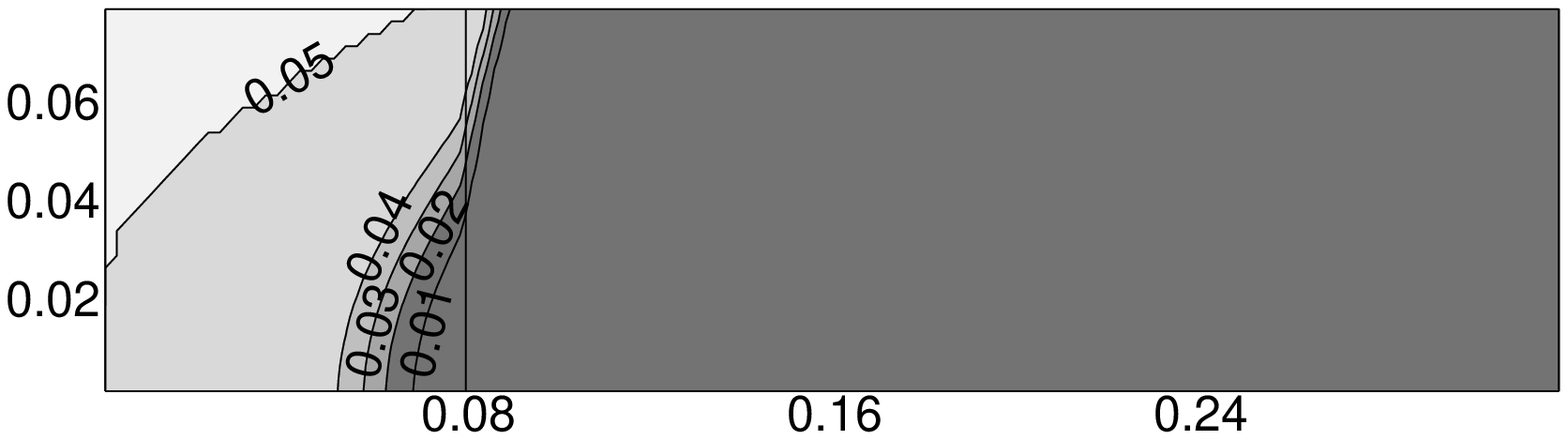,width=0.46\textwidth} \\[1ex]
\psfig{figure=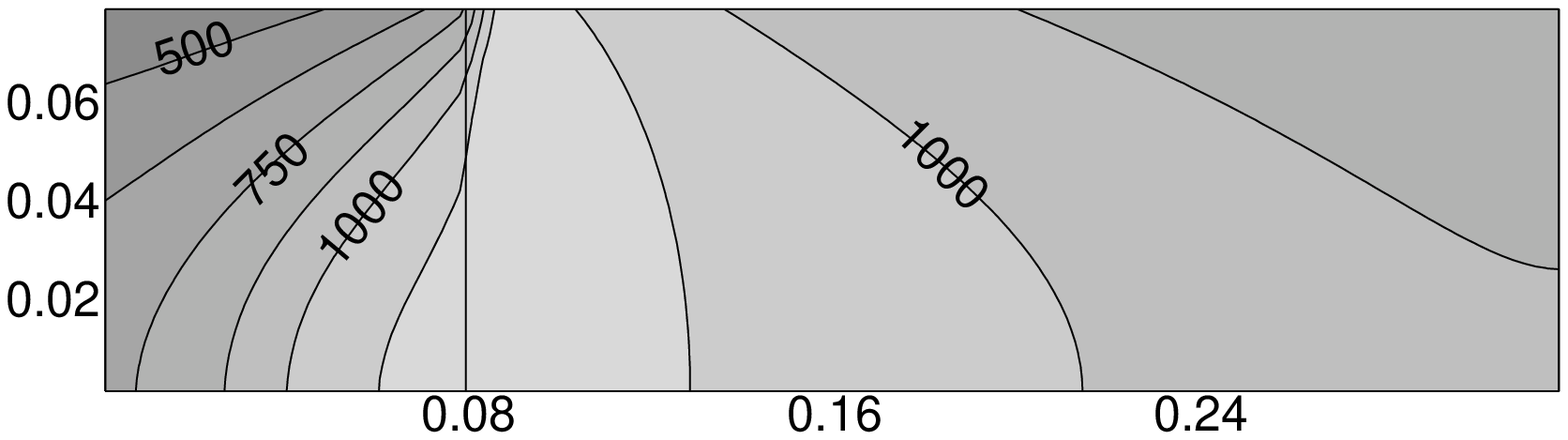,width=0.46\textwidth} &
\psfig{figure=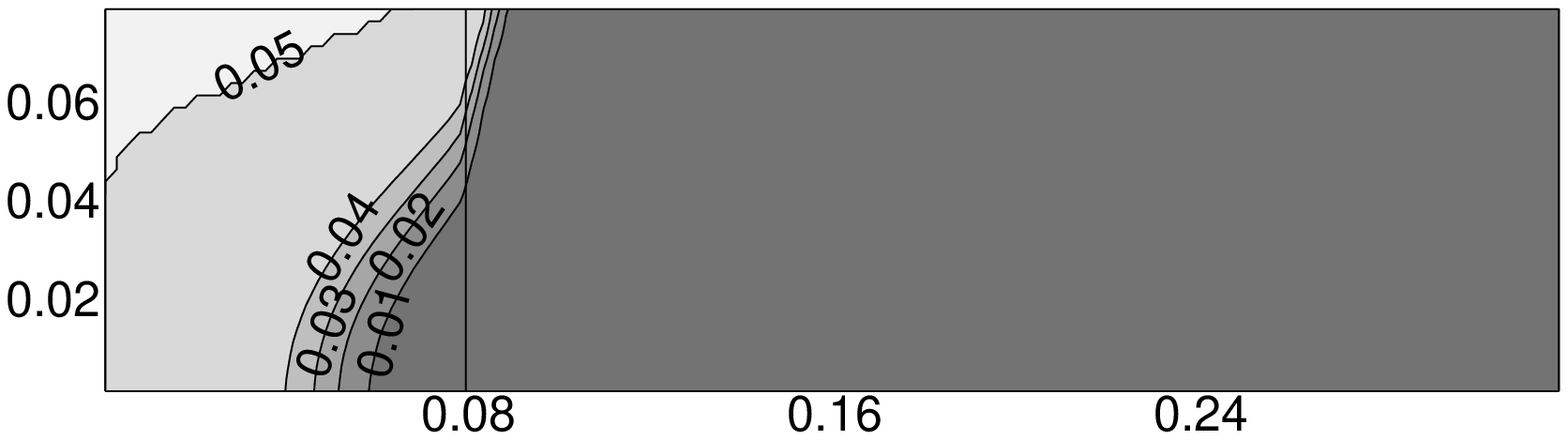,width=0.46\textwidth} \\[1ex]
(a) & (b)
\end{tabular}
\end{center}
\caption{Temperature (a) and mass fraction of reactant (b) for varying porosity}
\label{fig:sol:diff.por}
\end{figure}

\paragraph{a) Results for different values of $\phi_{\mathrm{l}}$ 
	      and $\phi_{\mathrm{h}}$ according to (\ref{diff.por})}
\hfill\\[1ex]
Figure \ref{fig:sol:diff.por} shows the calculated distributions of the temperature $T$ 
and the mass fraction of the reactant $y$ for
$t=150$, $t=300$, $t=500$, $t=1000$ and $t=5000$.

After switching off the external heat source,
the reaction zone, where the values of $y$ decrease from $0.05$ to $0.0$,
widens in direction of the symmetry line,
until nearly straight reaction front is generated.
In the following, the reaction zone moves towards the inflow boundary.
After reaching the region with lower porosity ($x_1 < 0.08$) this movement
slows down more and more.
Finally (for $t > 5000$) a stationary state is achieved, 
where the reaction zone is located around the interface 
between the two regions with higher and lower porosity.
This behavior can be explained in the following manner:
In the region with low porosity, the effective heat conductivity 
$\lambda_{\mathrm{eff}}$ is larger than in the region with high porosity.
Therefore the heat released by the reaction is conducted faster towards
the cooled burner wall. 
This effect is enforced by the fact that the heat transfer coefficient
at $\Gamma_{\mathrm{C}1}$ (no isolation) is much higher 
than at $\Gamma_{\mathrm{C}2}$.
Furthermore the reaction rate is proportional to the porosity $\phi$.
Thus there is a stronger heat production by chemical reaction in the region 
with high porosity.

\begin{figure}[ht]
\begin{center}
\begin{tabular}{cc}
\psfig{figure=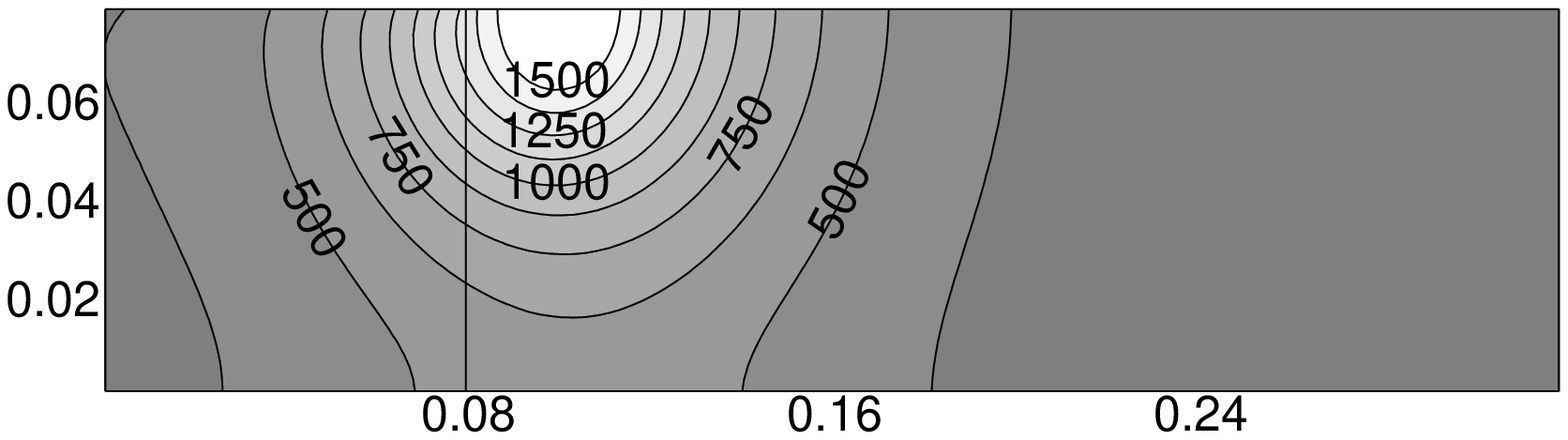,width=0.46\textwidth} &
\psfig{figure=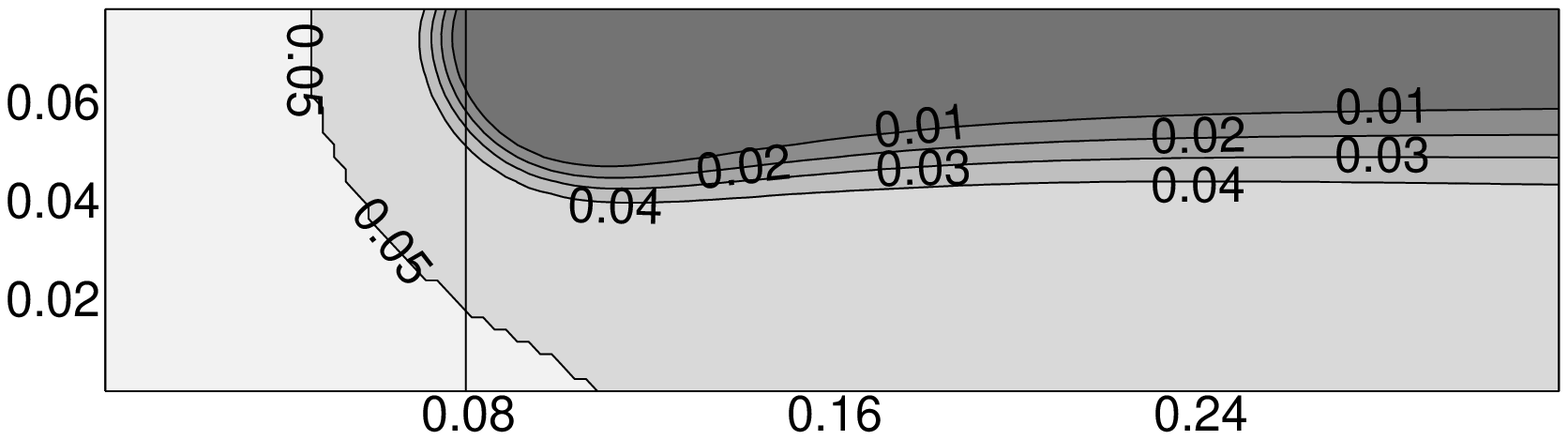,width=0.46\textwidth} \\[1ex]
\psfig{figure=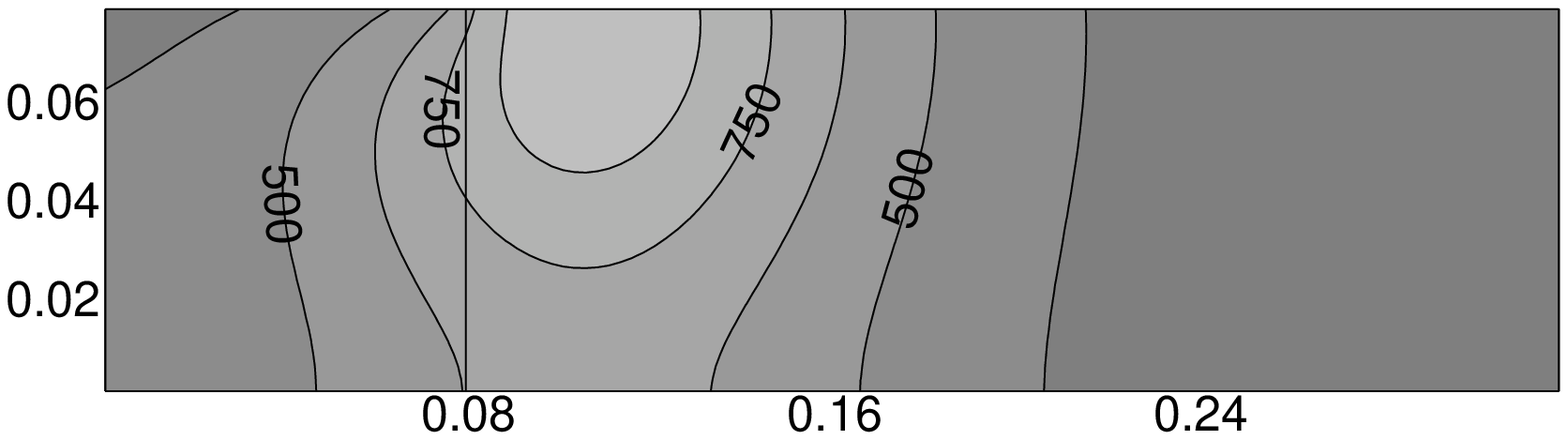,width=0.46\textwidth} &
\psfig{figure=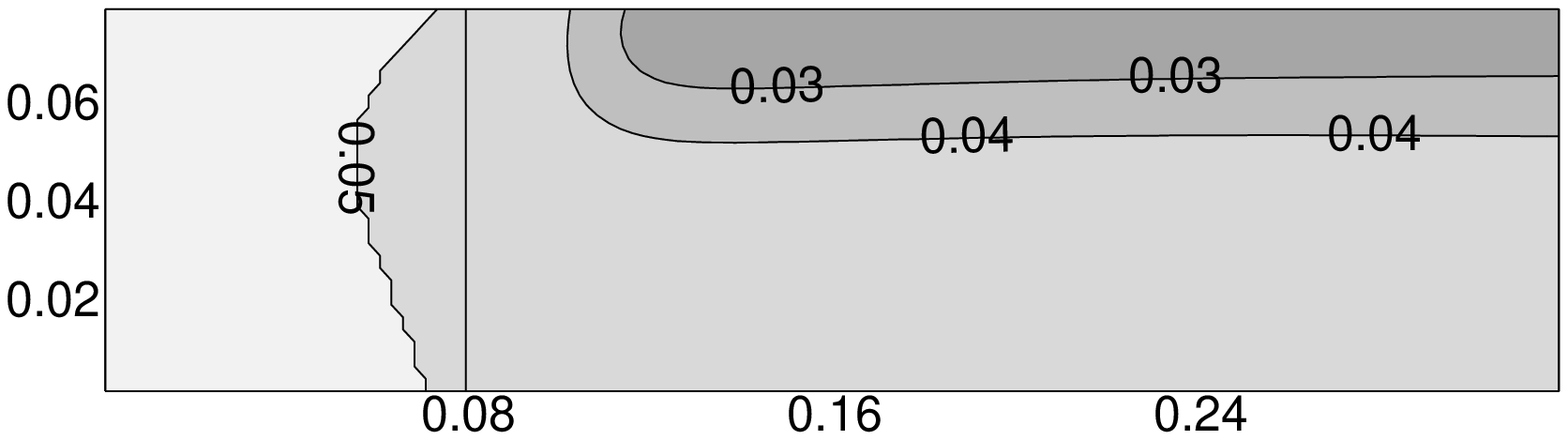,width=0.46\textwidth} \\[1ex]
\psfig{figure=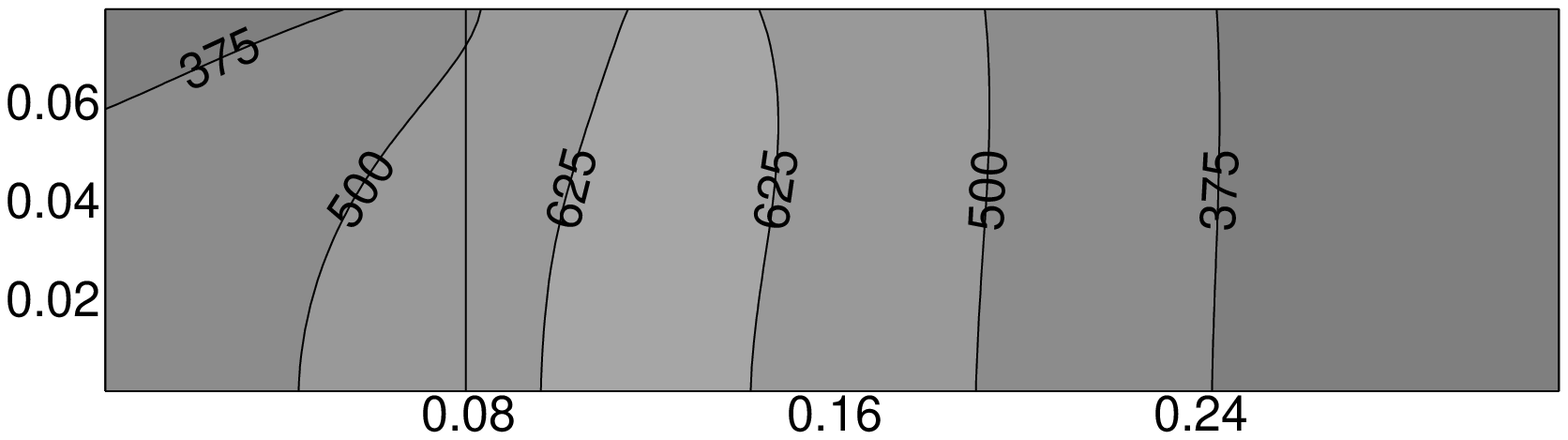,width=0.46\textwidth} &
\psfig{figure=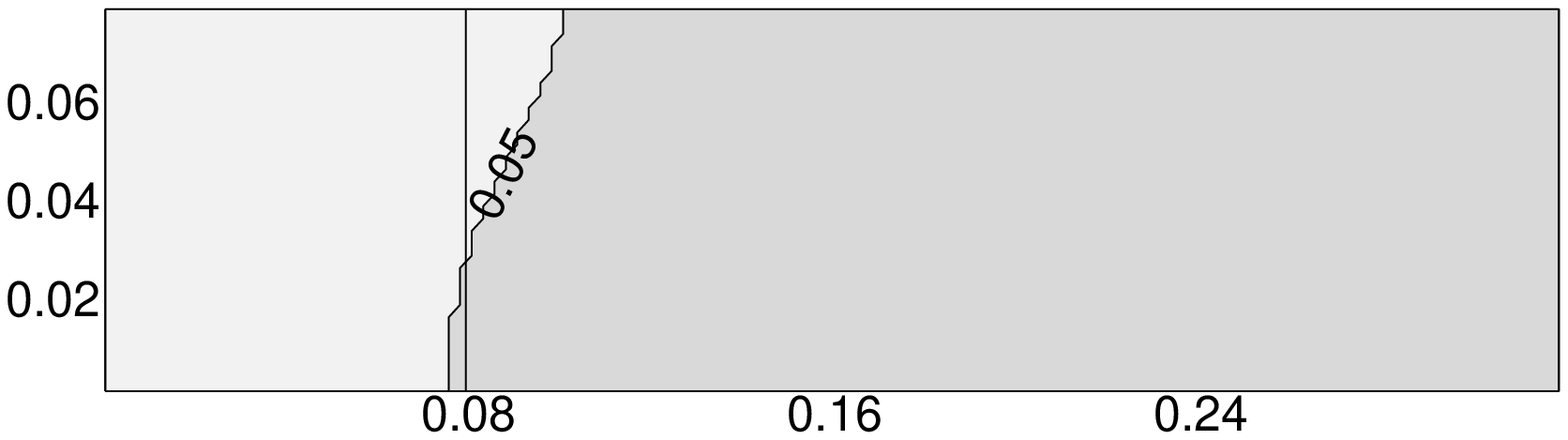,width=0.46\textwidth} \\[1ex]
\psfig{figure=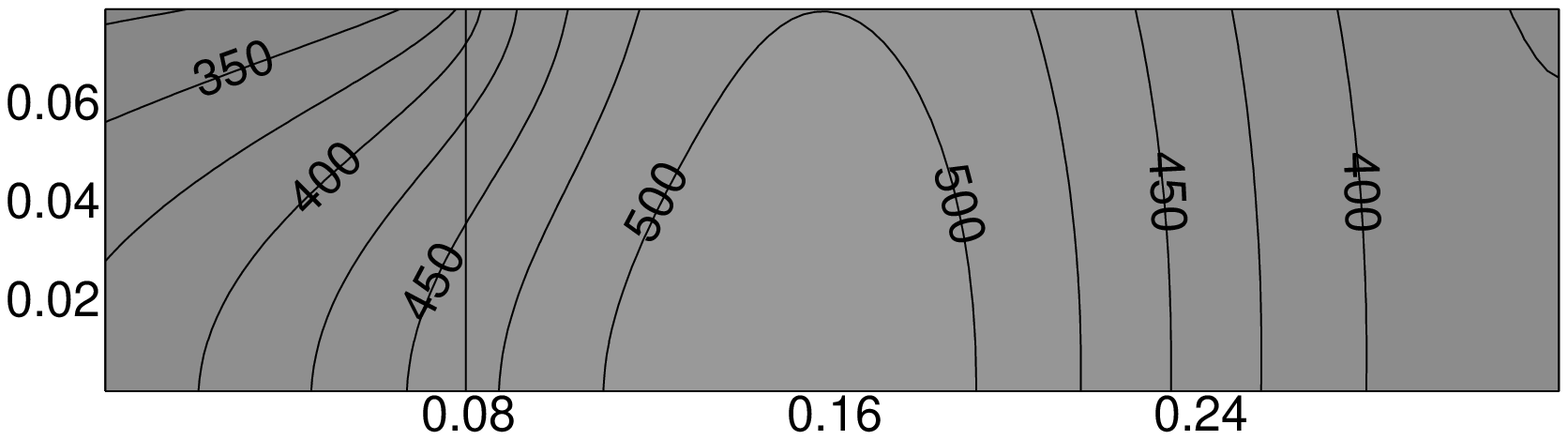,width=0.46\textwidth} &
\psfig{figure=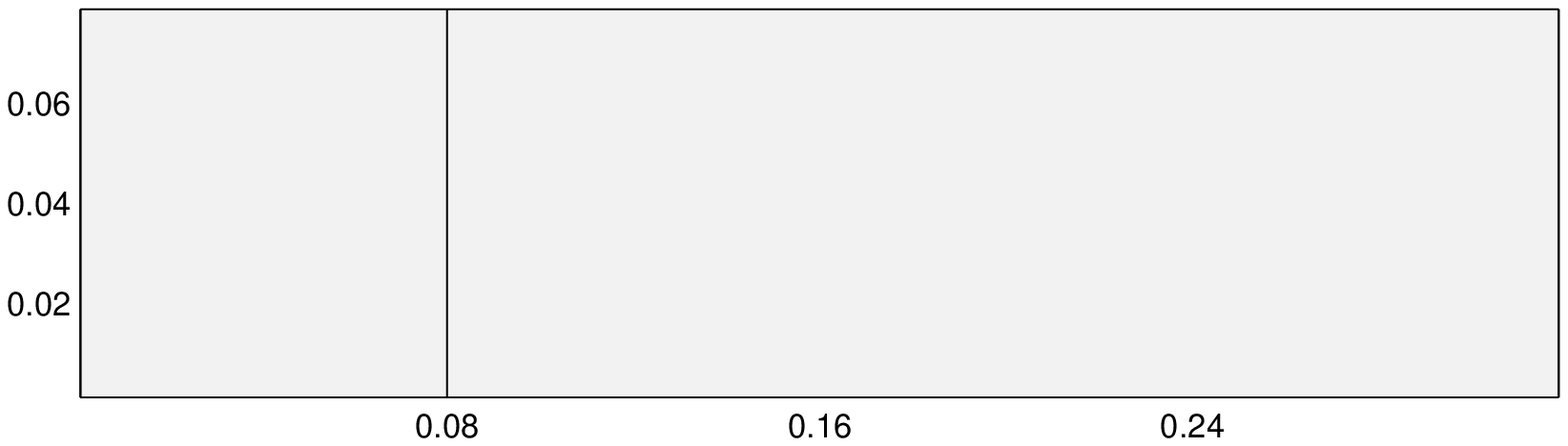,width=0.46\textwidth} \\[1ex]
\psfig{figure=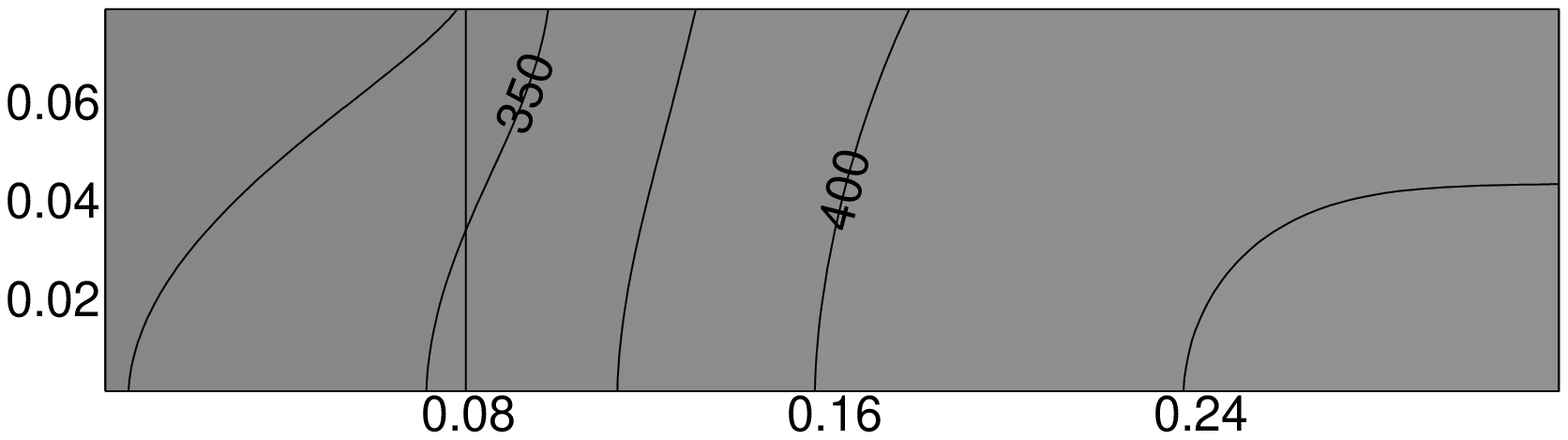,width=0.46\textwidth} &
\psfig{figure=fig/lp.spec500.ps,width=0.46\textwidth} \\[1ex]
(a) & (b)
\end{tabular}
\end{center}
\caption{Temperature (a) and mass fraction of reactant (b) for low porosity}
\label{fig:sol:low.por}
\end{figure}

\subparagraph{b) Results for fixed low value 
		 $\phi \equiv \phi_{\mathrm{l}}$ according to (\ref{low.por})}
\hfill\\[1ex]
Figures \ref{fig:sol:low.por} shows the calculated distributions of the temperature $T$ 
and the mass fraction of the reactant $y$ for
$t=150$, $t=200$, $t=300$, $t=500$ and $t=1000$.

After switching off the external heat source and starting the heat loss
across the burner wall, the reaction breaks down.
Again, this behavior is caused by the two effects mentioned in a) above.
On the one hand, there is less heat release of the reaction
due to lower porosity, 
on the other hand there is stronger heat conduction (towards the cooled wall)
due to the high portion of solid.
Hence the temperature falls quickly under the value needed for the reaction.
After the reaction is extinguished, there is only a flow problem coupled
with a heat conduction to solve. 
Finally the temperature $T$ is equal to the ambient temperature.

\begin{figure}[ht]
\begin{center}
\begin{tabular}{cc}
\psfig{figure=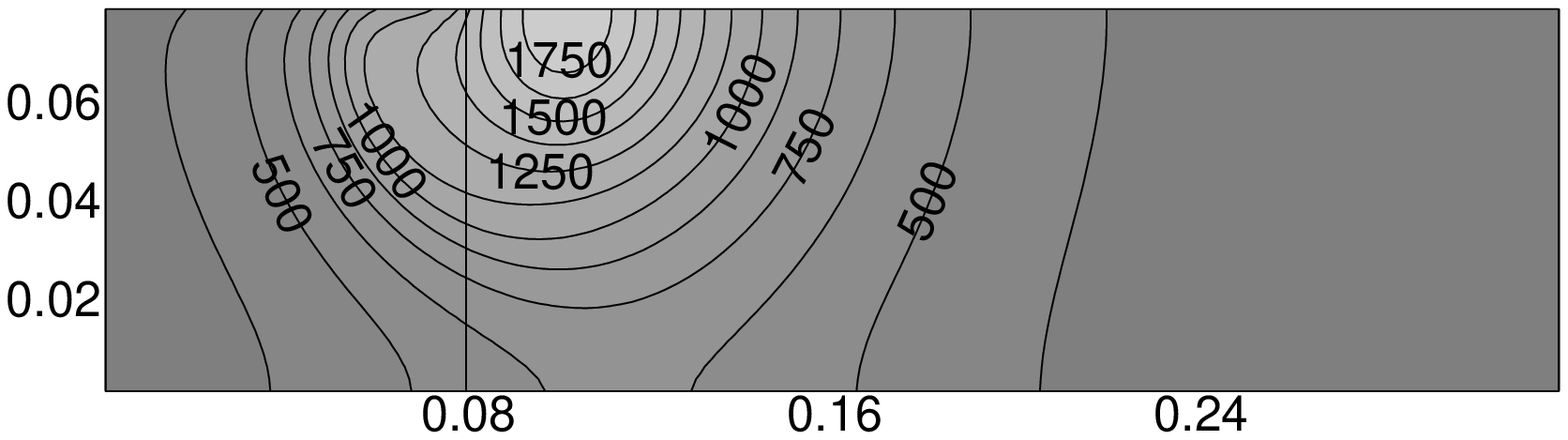,width=0.46\textwidth} &
\psfig{figure=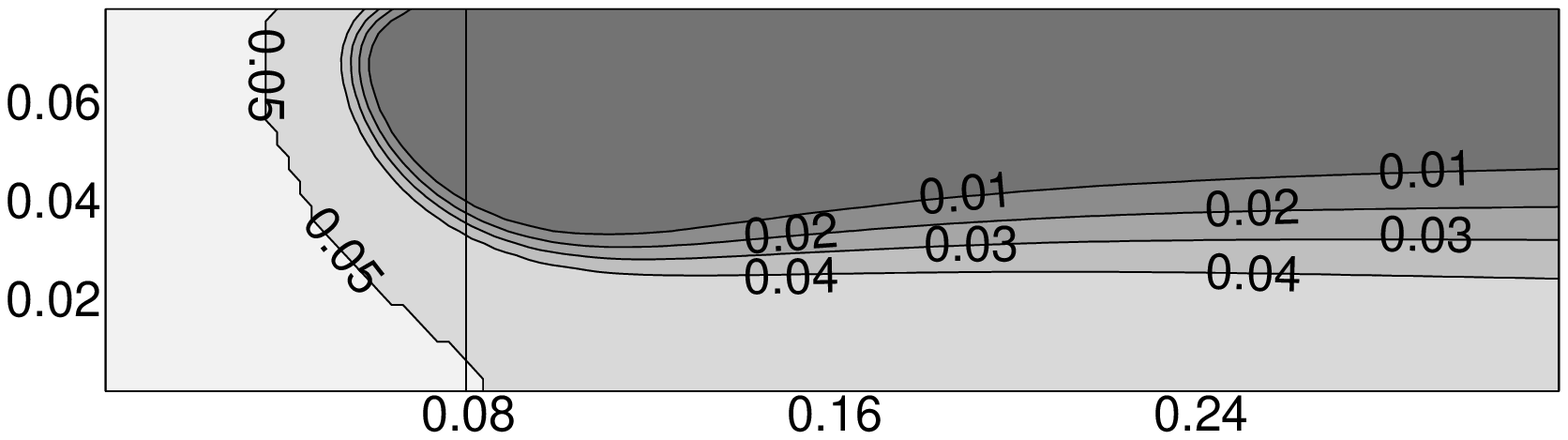,width=0.46\textwidth} \\[1ex]
\psfig{figure=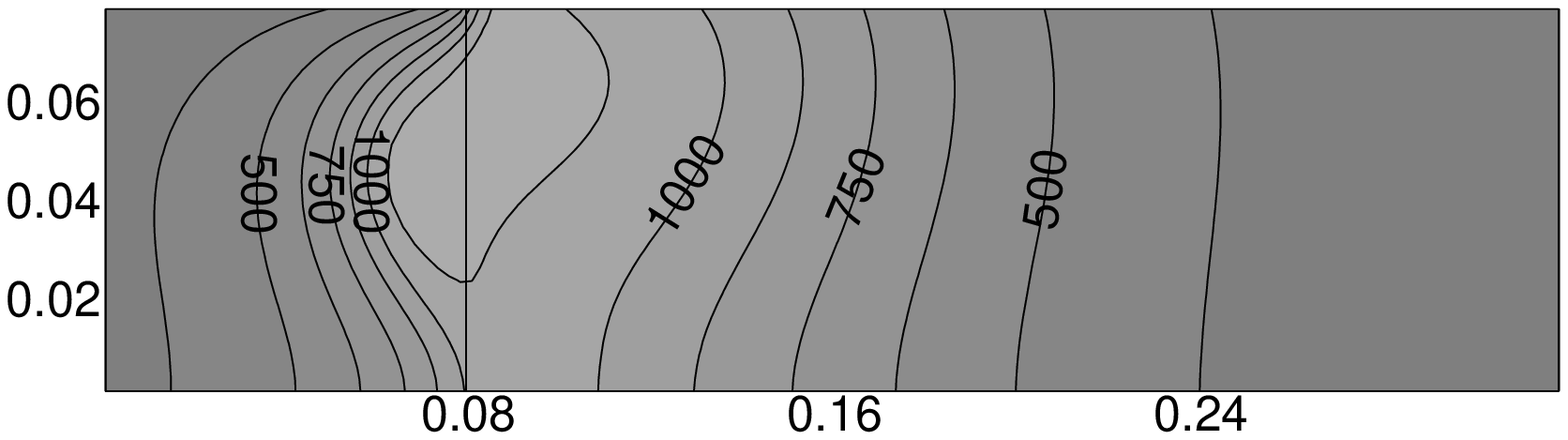,width=0.46\textwidth} &
\psfig{figure=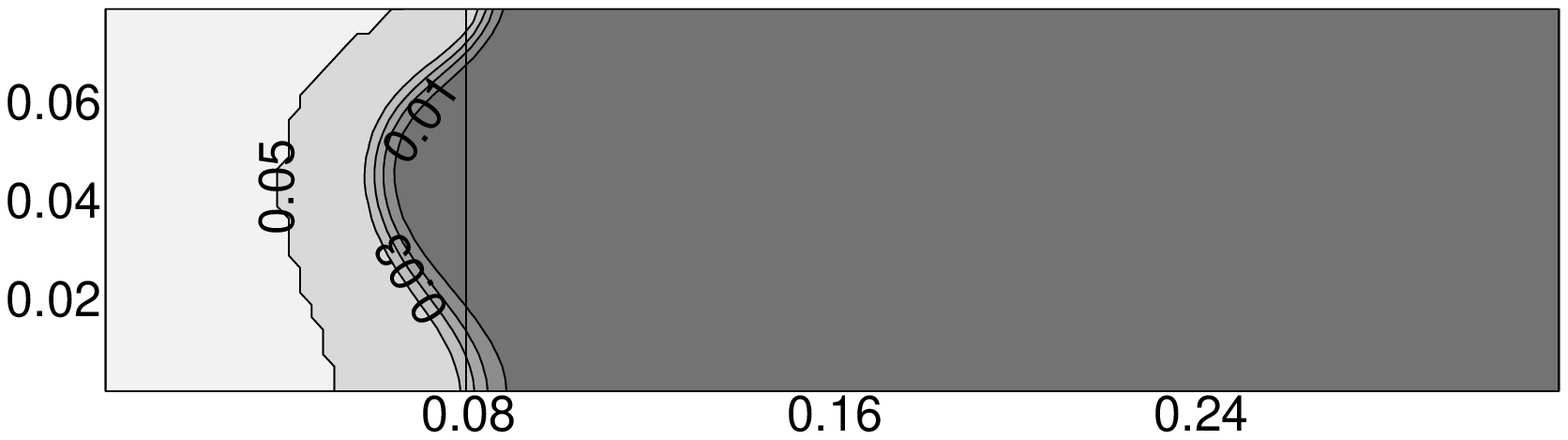,width=0.46\textwidth} \\[1ex]
\psfig{figure=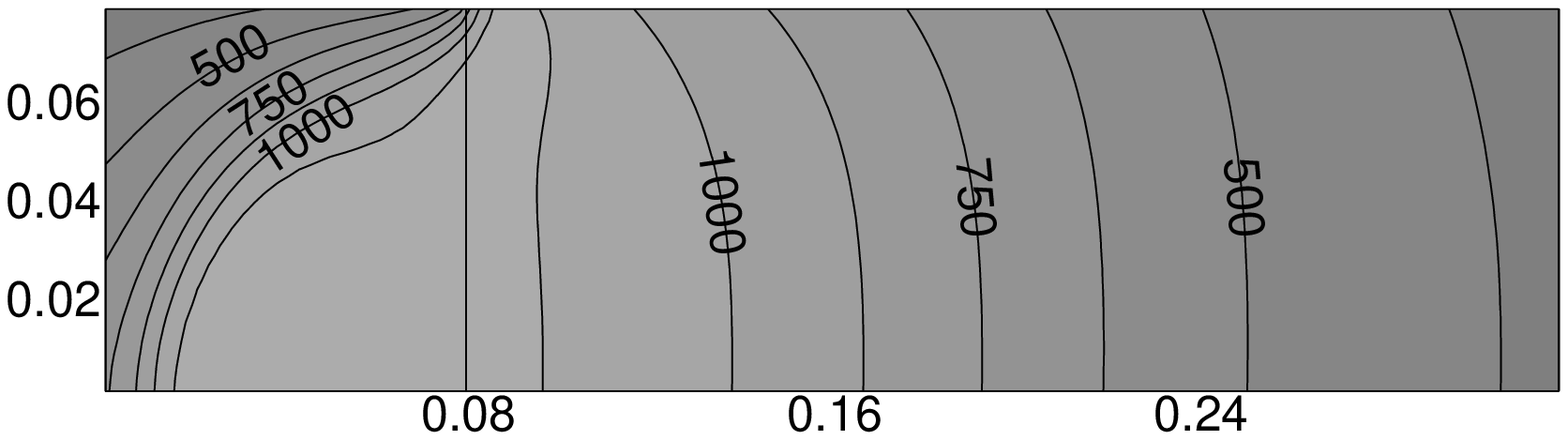,width=0.46\textwidth} &
\psfig{figure=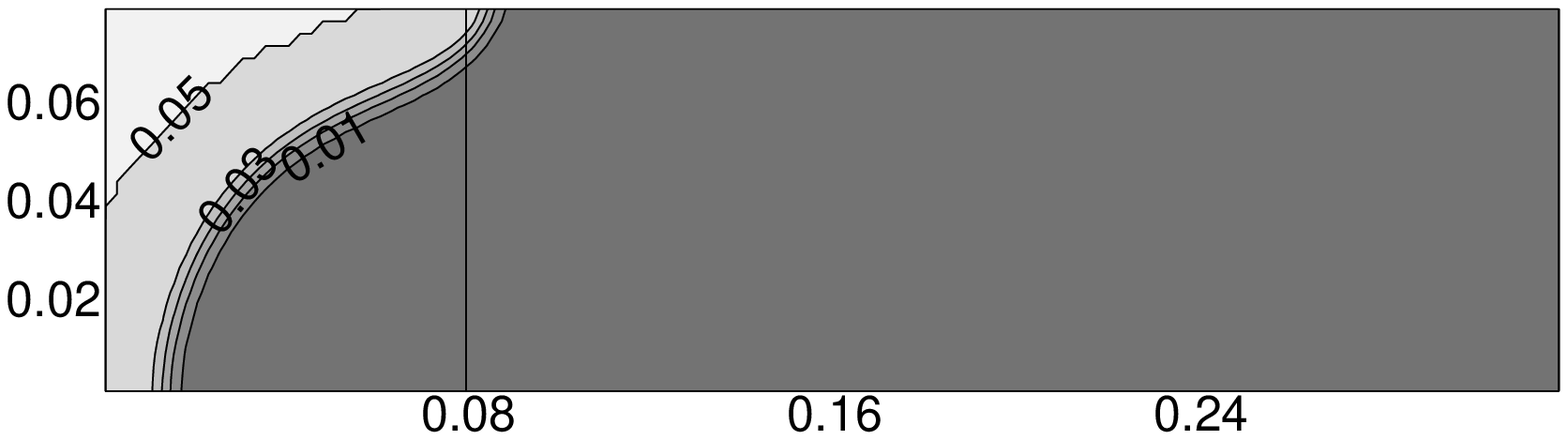,width=0.46\textwidth} \\[1ex]
\psfig{figure=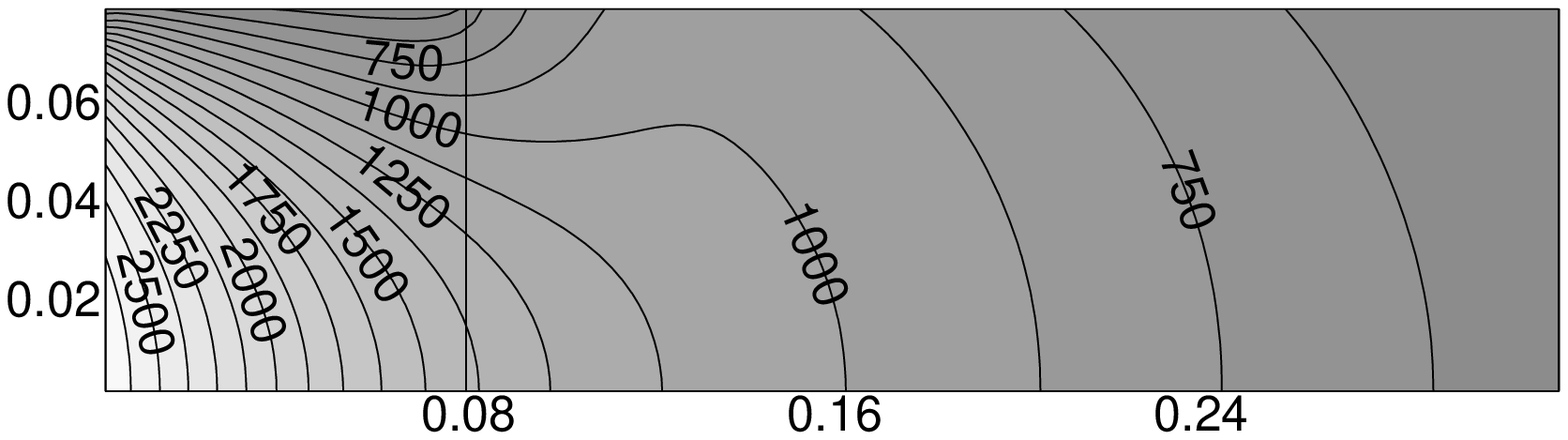,width=0.46\textwidth} &
\psfig{figure=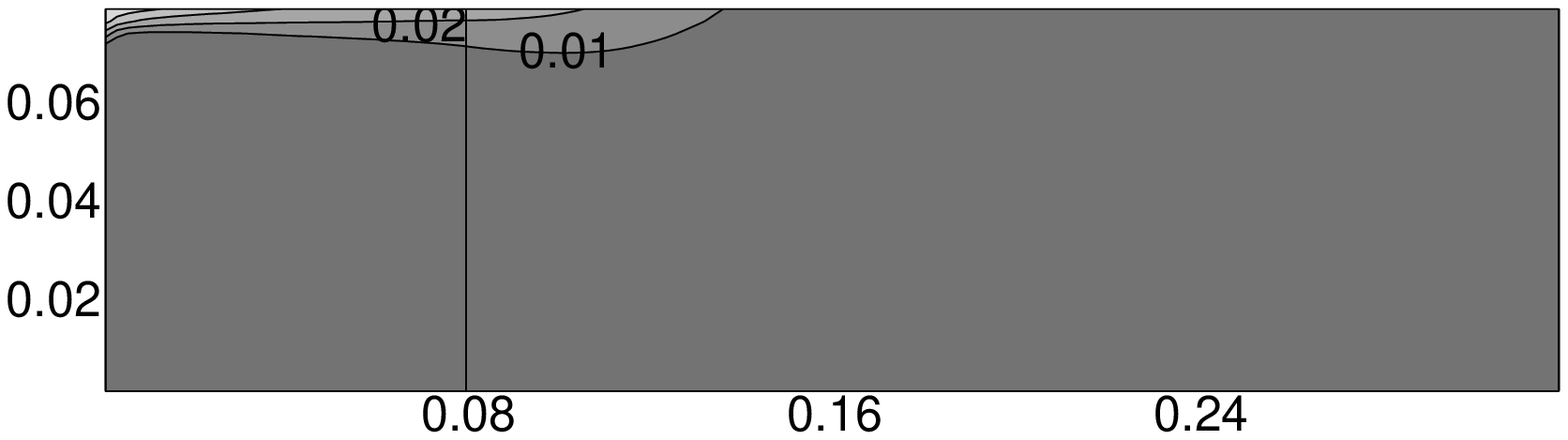,width=0.46\textwidth} \\[1ex]
\psfig{figure=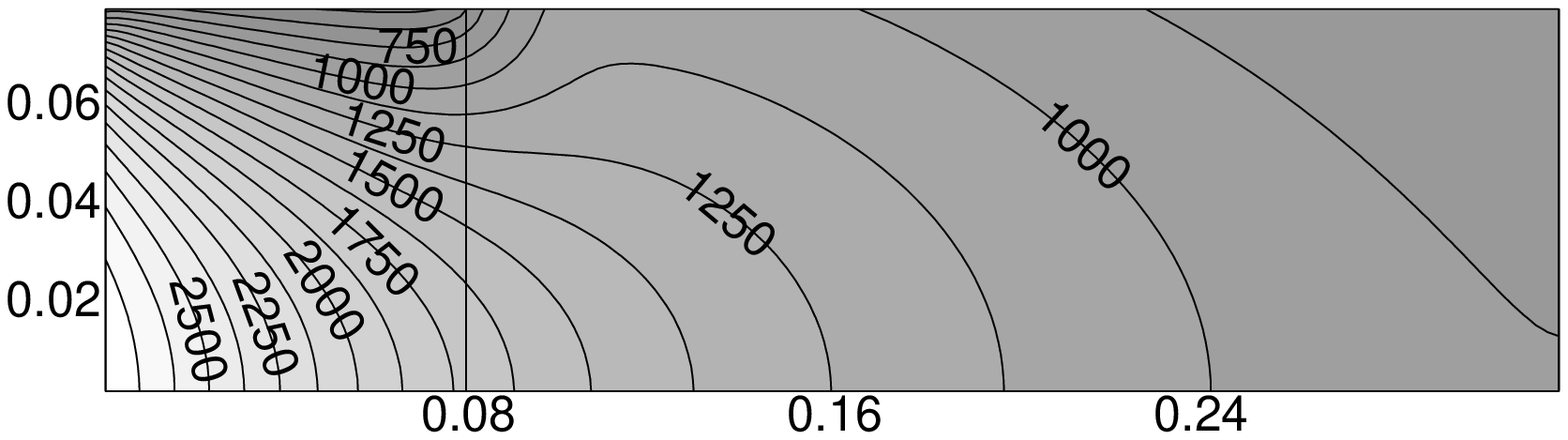,width=0.46\textwidth} &
\psfig{figure=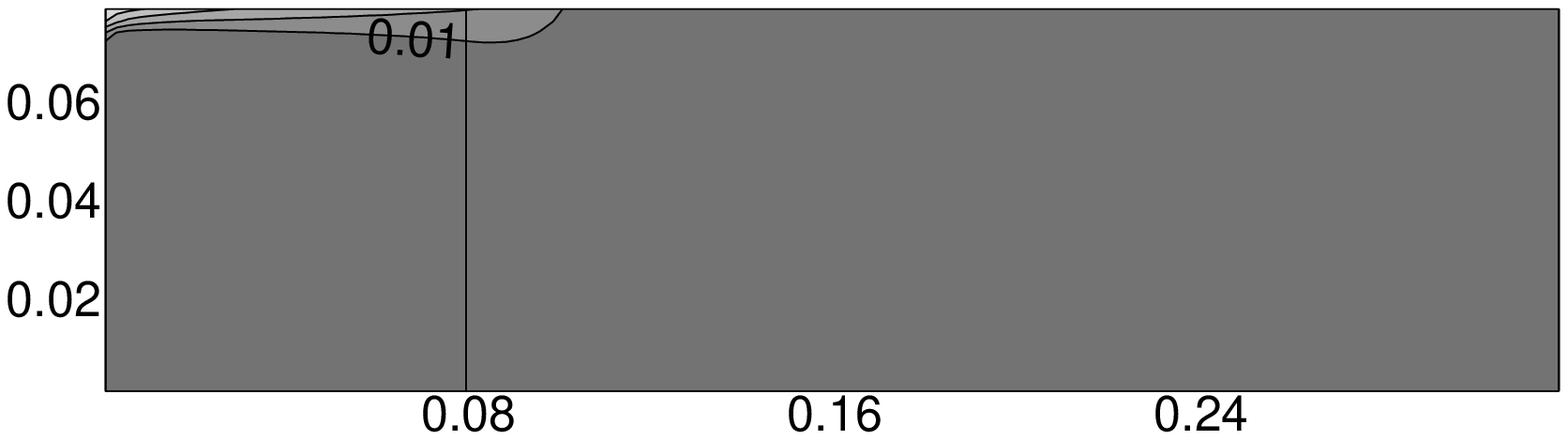,width=0.46\textwidth} \\[1ex]
(a) & (b)
\end{tabular}
\end{center}
\caption{Temperature (a) and mass fraction of reactant (b) for low porosity}
\label{fig:sol:high.por}
\end{figure}

\subparagraph{c) Results for fixed high value 
		 $\phi \equiv \phi_{\mathrm{h}}$ according to (\ref{high.por})}
\hfill\\[1ex]
Figure \ref{fig:sol:high.por} shows the calculated distributions of the temperature $T$ 
and the mass fraction of the reactant $y$ for
$t=100$, $t=200$, $t=300$, $t=500$ and $t=1000$.

In the beginning there is a behavior similar to case a).
But now, the intensity of the reaction in the preheating region
is not reduced due to lower porosity.
Likewise, the effective heat conductivity is not increased.
Hence the movement of the reaction zone towards the inflow boundary is not
retarded strong enough.
Finally, the reaction zone reaches the inflow boundary.

%
\section{Conclusion}
\label{sec:conc}
%
We proposed discretization and solution methods for the equations governing
combustion in porous inert media.
Decoupling the flow problem from the transport problems,
we can use different methods for the discretization of these two problems.
The discretization of the flow problem is performed 
by the mixed finite element method.
Thereby we obtain a good approximation for the mass flux, 
which governs the convective transport in the transport equations.
The transport problems are discretized by a cell-centered finite volume method.
The resulting nonlinear systems of equations for both problems are lineararized
with Newton's method, 
the linearized systems are solved with a multigrid algorithm.
Finally both subsystems are recoupled again in a Picard iteration.
Although this approach is presented on the basis of a strongly 
simplified model, it can be extended easily to more realistic situations.
Numerical simulations with this simplified model showed
that the stabilization of the reaction zone inside the porous burner
can be explained using only macroscopic quantities like porosity.

%
\section*{References}
%

\end{document}